\documentclass[11pt]{article}
\usepackage{amsmath}
\usepackage{amsfonts}
\usepackage{graphicx}
\usepackage{amssymb}
\usepackage{leftidx}
\usepackage{color}

\newtheorem{condition**}{A*}
\newtheorem{condition***}{C*}
\newtheorem{condition*}{C}

\newtheorem{proposition}{Proposition}[section]

\newtheorem{definition}{Definition}[section]
\newtheorem{theorem}{Theorem}[section]
\newtheorem{lemma}{Lemma}[section]
\newtheorem{remark}{Remark}[section]

\begin{document}

\title{A BSDE Approach to Stochastic Differential Games Involving Impulse
Controls and HJBI Equation}
\author{Liangquan Zhang$^{1}$\thanks{%
L. Zhang acknowledges the financial support partly by the National Nature
Science Foundation of China(Grant No. 11701040, 11871010 \&61871058) and the
Fundamental Research Funds for the Central Universities (No. 2019XD-A11).
E-mail: xiaoquan51011@163.com.} \\
{\small 1. School of Science} \\
{\small Beijing University of Posts and Telecommunications} \\
{\small Beijing 100876, China} }
\maketitle

\begin{abstract}
This paper focuses on zero-sum stochastic differential games in the
framework of forward-backward stochastic differential equations on a finite
time horizon with both players adopting impulse controls. By means of BSDE
methods, in particular that of the notion from Peng's stochastic \textit{%
backward semigroups}, we prove a dynamic programming principle for both the
upper and the lower value functions of the game. The upper and the lower
value functions are then shown to be the unique viscosity solutions of the
Hamilton-Jacobi-Bellman-Isaacs equations with a double-obstacle. As a
consequence, the uniqueness implies that the upper and lower value functions
coincide and the game admits a value.
\end{abstract}



\section{Introduction}

\label{sect:1}

Fleming and Souganidis \cite{FS} first investigated two-player zero-sum
stochastic differential games as a pioneering work in a rigorous manner and
proved that the lower and the upper value functions of such games fulfill
the dynamic programming principle shown to be the unique viscosity solutions
of the associated HJBI equations and coincide under the Isaacs condition.
This work developed the former results on differential games by Isaacs \cite%
{Isa}, Elliott and Kalton \cite{EK}, Friedman \cite{Fri}, Evans and
Souganidis \cite{ES} from the purely deterministic into the stochastic
framework and has made a huge progress in the field of stochastic
differential games. There are many works which extend the Fleming and
Souganidis approach into new contexts. For instance, Buckdahn, Cardaliaguet
and Rainer \cite{BCR} prove the existence of Nash equilibrium points for
stochastic nonzero-sum differential games and characterize them. Meanwhile,
the theory of backward stochastic differential equations (BSDE) has been
used to study stochastic differential games. In this direction, the reader
can see Hamad\`{e}ne and Lepeltier \cite{HL2} and Hamad\`{e}ne, Lepeltier,
and Peng \cite{HLP}. In particular, Buckdahn and Li \cite{BuckLi} developed
the findings obtained in \cite{HL2, HLP} and generalized the framework in
\cite{FS} to BSDE.
\textcolor[rgb]{1.00,0.00,0.00}{Another direction for generalization can be
seen in Bayraktar and Poor \cite{BP} and Browne \cite{Bro}. Concerning
optimal stopping, the interested reader is referred to the work of Ekstr\"{o}m and Peskir \cite{EP}, of Karatzas and Sudderth \cite{KS}, and of Karatzas
and Zamfirescu \cite{KZ}.}

Different from continuous control, impulse control is also an interesting
topic in stochastic control theory. There are three approaches to exploit
it. For functional analysis methods, see Bensoussan and Lions \cite{BL}. For
direct probabilistic methods, see Robin \cite{Robin} and Stettner \cite{Ste1}%
. For viscosity solution approaches concerning the study of impulse control,
see Lenhart \cite{Len}, Tang and Yong \cite{TY}, and Kharroubi et al. \cite%
{KMH}. Impulse control has been found as a useful tool for realistic models
in mathematical finance, for instance, transaction costs and liquidity risk.
For more information on this direction refer, in particular, Korn \cite{Korn}%
, Ly Vath, Mnif and Pham \cite{MP} and Bruder and Pham \cite{BP2}, Tang and
Hou \cite{TH}. For nonzero sum impulse control games, see \cite{ABCC}.
\textcolor[rgb]{1.00,0.00,0.00}{Some recent advances in numerical impulse
control can be seen in \cite{ARS,AF, ABL, Za}.} Besides, the celebrated Pontryagin's maximum principle for stochastic differential games within the framework of BSDE can be found in Wang and Yu \cite{WY2012} and Wang, Xiao and Xiong \cite{WXX2018}. The related topic in this fields see Chen and Wu \cite{CW2020} and Xu \cite{Xu2020}.

Cosso \cite{co} and El Asri and Mazid \cite{AM} studied a two-player
zero-sum stochastic differential game, with both players adopting impulse
controls on a finite time horizon of the following type: The state process
is governed by a $n$-dimensional SDE of the following type:%
\begin{eqnarray}
X_{s}^{t,x;u,v} &=&x+\int_{t}^{s}b\left( s,X_{s}^{t,x;u,v}\right) \mathrm{d}%
s+\int_{t}^{s}\sigma \left( s,X_{s}^{t,x;u,v}\right) \mathrm{d}%
W_{s}+\sum_{l\geq 1}\eta _{l}\mathbf{1}_{\left[ \rho _{l},T\right] }\left(
s\right)  \notag \\
&&+\sum_{m\geq 1}\xi _{m}\mathbf{1}_{\left[ \tau _{m},T\right] }\left(
s\right) \prod_{l\geq 1}\mathbf{1}_{\left\{ \tau _{m}\neq \rho _{l}\right\}
},  \label{SDE1}
\end{eqnarray}%
for all $s\in \left[ t,T\right] ,$ $P$-a.s., with $X_{t-}=x,$ on some
filtered probability space $\left( \Omega ,\mathcal{F},P\right) $, where $b:%
\left[ 0,T\right] \times \mathbb{R}^{n}\rightarrow \mathbb{R}^{n},$ $\sigma
\left( \cdot ,\cdot \right) :\left[ 0,T\right] \times \mathbb{R}%
^{n}\rightarrow \mathbb{R}^{n\times d}$ are given deterministic functions, $%
\left( W_{s}\right) _{s\geq 0}$ is an $d$-dimensional Brownian motion, $%
\left( x,t\right) $ are initial time and state. $u=\sum_{m\geq 1}\xi _{m}%
\mathbf{1}_{\left[ \tau _{m},T\right] }$ and $v=\sum_{l\geq 1}\eta _{l}%
\mathbf{1}_{\left[ \rho _{l},T\right] }$ are the impulse controls of player
I and player II, respectively. The infinite product $\Pi _{l\geq 1}\mathbf{1}%
_{\left\{ \tau _{m}\neq \rho _{l}\right\} }$ has the following meaning:
whenever the two players act together on the system at the same time, we
take into account only the action of player II.

The gain functional for player I (resp., cost functional for player II) of
the stochastic differential game is given by%
\begin{eqnarray}
J\left( t,x;u,v\right) &=&\mathbb{E}\Bigg [\int_{t}^{T}f\left(
s,X_{s}^{t,x;u,v}\right) \mathrm{d}s-\sum_{m\geq 1}c\left( \tau _{m},\xi
_{m}\right) \mathbf{1}_{\left\{ \tau _{m}\leq T\right\} }\left( s\right)
\prod_{l\geq 1}\mathbf{1}_{\left\{ \tau _{m}\neq \rho _{l}\right\} }  \notag
\\
&&+\sum_{l\geq 1}\chi \left( \rho _{l},\eta _{l}\right) \mathbf{1}_{\left\{
\rho _{l}\leq T\right\} }+g\left( X_{T}^{t,x;u,v}\right) \Bigg ],
\label{cost1}
\end{eqnarray}%
where $f:\left[ 0,T\right] \times \mathbb{R}^{n}\rightarrow \mathbb{R}$ and $%
g:\mathbb{R}^{n}\rightarrow \mathbb{R}$ are two given deterministic
functions, $f$ denoting the running function and $g$ the payoff. The
function $c$ is the cost function for player I and the gain function for
player II, representing that when player I performs an action he/she has to
pay a cost, resulting in a gain for player II. Analogously, $\chi $ is the
cost function for player II and the gain function for player I. Cosso \cite%
{co} and El Asri and Mazid \cite{AM} (under weak assumptions) have shown
that the upper and lower value functions coincide and the game admits a
value.

The theory of BSDE can be traced back to Bismut \cite{BJ} who studied linear
BSDE motivated by stochastic control problems. Pardoux and Peng \cite{PP1}
proved the well-posedness for nonlinear BSDE. Subsequently, Duffie and
Epstein \cite{DE} introduced the notion of recursive utilities in continuous
time, which is actually a type of BSDE where the generator $f$ is
independent of $z$. Then, El Karoui et al. \cite{KPQ} extended the recursive
utility to the case where $f$ contains $z$. The term $z$ can be interpreted
as an ambiguity aversion term in the market (see Chen and Epstein 2002 \cite%
{CZJ}). Particularly, the celebrated Black-Scholes formula indeed provided
an effective way of representing the option price (which is the solution to
a kind of linear BSDE) through the solution of the Black-Scholes equation.
Since then, BSDE have been extensively studied and used in the areas of
applied probability and optimal stochastic controls, particularly in
financial engineering (cf. \cite{KPQ}).

In our present work, employing BSDE methods, in particular, the notion of
stochastic backward semigroups (Peng \cite{P3}), allows us to prove the
dynamic programming principle for the upper and lower value functions of the
game, with both players adopting impulse controls on a finite time horizon,
and to derive from it with the help of Peng's method (similar to \cite{P3,
P2}) the associated HJBI equations with a double-obstacle. To the best of
our knowledge, this is the first work studying impulse control games via
BSDE.

Consider
\begin{eqnarray}
Y_{s}^{t,x;u,v} &=&\Phi \left( X_{T}^{t,x;u,v}\right) +\int_{s}^{T}f\left(
r,X_{r}^{t,x;u,v},Y_{r}^{t,x;u,v},Z_{r}^{t,x;u,v}\right) \mathrm{d}%
r-\int_{s}^{T}Z_{r}^{t,x;u,v}\mathrm{d}W_{r}  \notag \\
&&-\sum_{m\geq 1}c\left( \tau _{m},\xi _{m}\right) \mathbf{1}_{\left\{ \tau
_{m}\leq T\right\} }\prod_{l\geq 1}\mathbf{1}_{\left\{ \tau _{m}\neq \rho
_{l}\right\} }+\sum_{l\geq 1}\chi \left( \rho _{l},\eta _{l}\right) \mathbf{1%
}_{\left\{ \rho _{l}\leq T\right\} },  \label{BSDE1}
\end{eqnarray}%
where $X_{\cdot }^{t,x;u,v}$ is defined in (\ref{SDE1}). The existence and
uniqueness of BSDE (\ref{BSDE1}) under certain conditions can be guaranteed
in the next section. We are interested in studying two-player zero-sum
stochastic differential game, with both players adopting impulse controls on
a finite time horizon driven by FBSDEs (\ref{SDE1})-(\ref{BSDE1}). Compared with above
literature, our paper has several new features. The novelty of the
formulation and the contribution in this paper can be stated as follows:

\begin{itemize}
\item First, in the framework of BSDE, the terminal condition will turn out
to be a traditional $\Phi (X_{T}^{t,x;u,v})$ plus gains functions with
impulses controls. This new trait makes the backward semigroup valid and
avoids the It\^{o}'s formula with jumps.

\item Second, in Cosso \cite{co} and El Asri and Mazid \cite{AM}, the cost
functional is defined by (\ref{cost1}) via linear expectation. Our paper
considers a more general running cost functional, which implies that the
cost functionals will be supported by a BSDE, which in fact defines a
nonlinear expectation.

\item At last, as a response to one of Cosso's closing comments (Cosso \cite%
{co}, 2013): \textquotedblleft \textit{we could apply backward stochastic
differential equations methods to provide a probabilistic representation,
known as the nonlinear Feynman--Kac formula, for the value function of the
game\textquotedblright }, our paper aims to further perfect this theory in
the framework of BSDE, with more rigorous proofs and new techniques
introduced.
\end{itemize}

The rest of this paper is organized as follows: after some preliminaries and
notations in the second section, we devote the third section to studying the
regularity properties of the upper and lower value functions. Moreover, we
prove the dynamic programming principle for the stochastic differential game
with some corollaries and generalizations, which are useful in proving that
the two value functions are viscosity solutions to the HJBI equation in
Section \ref{sec4}. Furthermore, under certain assumptions, we establish the
comparison theorem for the HJBI equation, from which one may deduce that the
game admits a value. Finally, in Section \ref{sec5}, we conclude and
schedule possible generalizations in future. Some proofs can be found in
Appendix \ref{APP}.

\section{Preliminaries and Notations}

\label{sect2}

Throughout this paper, we denote by $\mathbb{R}^{n}$ the space of $n$%
-dimensional Euclidean space, by $\mathbb{R}^{n\times d}$ the space the
matrices with order $n\times d$. The probability space is the classical
Wiener space $\left( \Omega ,\mathcal{F},P\right) $, and the Brownian motion
$W$ will be the coordinate process on $\Omega $. Precisely: $\Omega $ is the
set of continuous functions from $\left[ 0,T\right] $ to $\mathbb{R}^{d}$
starting from $0$ ($\Omega =C_{0}\left( \left[ 0,T\right] ;\mathbb{R}%
^{d}\right) $), $\mathcal{F}$ is the Borel $\sigma $-algebra over $\Omega $,
completed with respect to the Wiener measure $P$ on this space, and $W$
denotes the coordinate process: $W_{s}\left( \omega \right) =\omega _{s}$, $%
s\in \left[ 0,T\right] $, $\omega \in \Omega $. By $\mathcal{F=}\left\{
\mathcal{F}_{s},0\leq s\leq T\right\} $, we denote the natural filtration
generated by $\left\{ W_{s}\right\} _{0\leq s\leq T}$ and augmented by all
P-null sets, i.e., $\mathcal{F}_{s}=\sigma \left\{ W_{r},r\leq s\right\}
\vee \mathcal{N}_{P},$ $s\in \left[ t,T\right] ,$ where $\mathcal{N}_{P}$ is
the set of all $P$-null subsets and $T>0$ a fixed real time horizon. For
each $t>0$, we denote by $\left\{ \mathcal{F}_{s}^{t},t\leq s\leq T\right\} $
the natural filtration of the Brownian motion $\left\{ W_{s}-W_{t},\text{ }%
t\leq s\leq T\right\} $, augmented by $\mathcal{N}_{P}$. $\top $ appearing
in this paper as superscript denotes the transpose of a matrix. $U$ and $V$
are two convex cones of $\mathbb{R}^{n}$ with $U\subset V$. In what follows,
$C$ represents a generic constant, which can be different from line to line.

Now, we give the following definition.

\begin{definition}
\label{s1}An \emph{impulse control} $u=\sum_{m\geq 1}\xi _{m}\mathbf{1}_{%
\left[ \tau _{m},T\right] }$ for player I (resp. $v=\sum_{l\geq 1}\eta _{l}%
\mathbf{1}_{\left[ \rho _{l},T\right] }$ for player II) on $\left[ t,T\right]
$ is such that

(1) $\left( \tau _{m}\right) _{m}$ (resp., $\left( \rho _{l}\right) _{l}$),
the action time, is a nondecreasing sequence of $\mathcal{F}$-stopping time,
valued in $\left[ t,T\right] \cup \left\{ +\infty \right\} .$

(2) $\left( \xi _{m}\right) _{m}$ (resp., $\left( \eta _{l}\right) _{l}$),
the actions, is a sequence of $U$-valued (resp., $V$-valued) random
variable, where each $\xi _{m}$ (resp., $\eta _{l}$) is $\mathcal{F}_{\tau
_{m}}$-measurable (resp., $\mathcal{F}_{\rho _{l}}$-measurable).
\end{definition}

\begin{remark}
Let $D\left( \left[ 0,T\right] ;\mathbb{R}^{m}\right) $ be the space of all
functions $\xi :\left[ 0,T\right] \rightarrow \mathbb{R}^{m}$ that are right
limit with left continuous. Then, the pure jump part of $\xi $ is defined by
$\xi ^{j}\left( t\right) =\sum_{0\leq s\leq t}\Delta \xi \left( s\right) ,$
and the continuous part is $\xi ^{c}\left( t\right) =\xi \left( t\right)
-\xi ^{j}\left( t\right) .$ By Lebesgue decomposition Theorem that we have $%
\xi ^{c}\left( t\right) =\xi ^{ac}\left( t\right) +\xi ^{sc}\left( t\right) $%
, $t\in \left[ 0,T\right] $, where $\xi ^{ac}\left( t\right) $ is called the
absolutely continuous part of $\xi ,$ and $\xi ^{sc}$ the singularly
continuous part of $\xi $. Thus, we obtain that%
\begin{equation*}
\xi \left( t\right) =\xi ^{ac}\left( t\right) +\xi ^{sc}\left( t\right) +\xi
^{j}\left( t\right) ,\text{ }t\in \left[ 0,T\right] ,\text{ unique!}
\end{equation*}%
If we assume that $\xi ^{ac}\left( t\right) +\xi ^{sc}\left( t\right) \equiv
0,$ $t\in \left[ 0,T\right] ,$ then the singular control performs a special
form of a pure jump process, so-called impulse control (see \cite{YZ} for
details).
\end{remark}

We now introduce the following spaces of processes:
\begin{align*}
\mathcal{S}^{2}(0,T;\mathbb{R})\triangleq & \left\{ \mathbb{R}^{n}\text{%
-valued }\mathcal{F}_{t}\text{-adapted process }\phi (t)\text{; }\mathbb{E}%
\left[ \sup\limits_{0\leq t\leq T}\left\vert \phi _{t}\right\vert ^{2}\right]
<\infty \right\} , \\
\mathcal{M}^{2}(0,T;\mathbb{R})\triangleq & \left\{ \mathbb{R}^{n}\text{%
-valued }\mathcal{F}_{t}\text{-adapted process }\varphi (t)\text{; }\mathbb{E%
}\left[ \int_{0}^{T}\left\vert \varphi _{t}\right\vert ^{2}\mbox{\rm d}t%
\right] <\infty \right\} ,
\end{align*}%
and denote $\mathcal{N}^{2}\left[ 0,T\right] =\mathcal{S}^{2}(0,T;\mathbb{R}%
^{n})\times \mathcal{S}^{2}(0,T;\mathbb{R})\times \mathcal{M}^{2}(0,T;%
\mathbb{R}^{n}).$ Clearly, $\mathcal{N}^{2}\left[ 0,T\right] $ forms a
Banach space.

We assume that the following conditions hold.

\begin{description}
\item[(A1)] The coefficients $b:[0,T]\times \mathbb{R}^{n}\rightarrow
\mathbb{R}^{n}$ and $\sigma :[0,T]\times \mathbb{R}^{n}\rightarrow \mathbb{R}%
^{n}$ are continuous on $[0,T]\times \mathbb{R}^{n},$ Lipschitz continuous
in the state variable $x,$ uniformly with respect to time, and bound on $%
[0,T]\times \mathbb{R}^{n}.$

\item[(A2)] The coefficients $f:[0,T]\times \mathbb{R}^{n}\times \mathbb{R}%
\times \mathbb{\mathbb{R}}^{d}$ and $\Phi :\mathbb{R}^{n}\rightarrow \mathbb{%
R}$ are continuous on $[0,T]\times \mathbb{R}^{n}\times \mathbb{R}\times
\mathbb{\mathbb{R}}^{d},$ Lipschitz continuous in the state variable $\left(
x,y,z\right) ,$ uniformly with respect to time, and bounded on $[0,T]\times
\mathbb{R}^{n}.$
\end{description}

To get a \emph{well-defined} gain functional, we add the following
assumption and introduce the concept of admissible impulse controls.
Meanwhile, to ensure that multiple impulses occurring at the same time are
suboptimal, we put (like in Cosso \cite{co}) the following:

\begin{description}
\item[(A3)] Let cost functions $c:\left[ 0,T\right] \times U\rightarrow
\mathbb{R}$ and $\chi :\left[ 0,T\right] \times V\rightarrow \mathbb{R}$ be
measurable and $1/2$ H\"{o}lder continuous in time, uniformly with respect
to the other variable. Furthermore,
\begin{equation}
\inf_{\left[ 0,T\right] \times U}c\left( t,\xi \right) >0,\text{ }\inf_{%
\left[ 0,T\right] \times V}\chi \left( t,\eta \right) >0,  \label{a1}
\end{equation}%
and there exists a function $h:\left[ 0,T\right] \rightarrow \left(
0,+\infty \right) $ such that for all $t\in \left[ 0,T\right] $,%
\begin{equation}
c\left( t,y_{1}+z+y_{2}\right) \leq c\left( t,y_{1}\right) -\chi \left(
t,z\right) +c\left( t,y_{2}\right) -h\left( t\right)  \label{a2}
\end{equation}%
and
\begin{equation}
\chi \left( t,z_{1}+z_{2}\right) \leq \chi \left( t,z_{1}\right) +\chi
\left( t,z_{2}\right) -h\left( t\right)  \label{a3}
\end{equation}%
for $y_{1},z,y_{2}\in U$ and $z_{1},z_{2}\in V.$ Moreover,
\begin{equation}
c\left( t,y\right) \geq c\left( \check{t},y\right) \text{ and }\chi \left(
t,y\right) \geq \chi \left( \check{t},y\right)  \label{a4}
\end{equation}%
for all $t,\check{t}\in \left[ 0,T\right] $ satisfying $t\leq \check{t},$ $%
y\in U$ and $z\in V.$
\end{description}

\begin{definition}
\label{s2}Let $u=\sum_{m\geq 1}\xi _{m}\mathbf{1}_{\left[ \tau _{m},T\right]
}$ be an impulse control on $\left[ t,T\right] $, and let $\sigma ,$ $\tau $
be two $\left[ t,T\right] $-valued $\mathcal{F}$-stopping times. Then we
define the \emph{restriction} $u_{\left[ \tau ,\sigma \right] }$ of the
impulse control $u$ by%
\begin{equation}
u_{\left[ \tau ,\sigma \right] }\left( s\right) =\sum_{m\geq 1}\xi _{\mu
_{t,\tau }\left( u\right) +m}\mathbf{1}_{\left\{ \tau _{\mu _{t,\tau }\left(
u\right) +m\leq s\leq \sigma }\right\} }\left( s\right) ,\text{ }\tau \leq
s\leq \sigma ,  \label{fimp}
\end{equation}%
where $\mu _{t,\tau }\left( u\right) $ is called the number of impulses up
to time $\tau ,$ namely, $\mu _{t,\tau }\left( u\right) :=\sum_{m\geq 1}%
\mathbf{1}_{\left\{ \tau _{m}\leq \tau \right\} }.$
\end{definition}

We now introduce the following subspaces of admissible controls.

\begin{definition}
\label{s3}An admissible impulse control $u$ for player I (resp., $v$ for
player II) on $\left[ t,T\right] $ is an impulse control for player I
(resp., II) on $\left[ t,T\right] $ with a finite average number of
impulses, i.e., $\mathbb{E}\left[ \mu _{t,T}\left( u\right) \right] <\infty $%
, resp., $\mathbb{E}\left[ \mu _{t,T}\left( v\right) \right] <\infty $, in
which $\mu _{t,T}\left( u\right) $ is given by (\ref{fimp}). The set of all
admissible impulse controls for player I (resp., II) on $\left[ t,T\right] $
is denoted by $\mathcal{U}_{t,T}$ (resp., $\mathcal{V}_{t,T}$). We identify
two impulse controls $u=\sum_{m\geq 1}\xi _{m}\mathbf{1}_{\left[ \tau _{m},T%
\right] }$ and $\tilde{u}=\sum_{m\geq 1}\tilde{\xi}_{m}\mathbf{1}_{\left[
\tilde{\tau}_{m},T\right] }$ in $\mathcal{U}_{t,T}$, and we write $u\equiv
\tilde{u}$ on $\left[ t,T\right] $ if $P\left( \left\{ u=\tilde{u},\text{
a.e. on }\left[ t,T\right] \right\} \right) =1$. Similarly, we interpret $%
u\equiv \tilde{u}$ on $\left[ t,T\right] $ in $\mathcal{V}_{t,T}$.
\end{definition}

Finally, we have still to define the admissible strategies for the game.

\begin{definition}
\label{s4}A nonanticipative strategy for player I on $\left[ t,T\right] $ is
a mapping $\alpha :\mathcal{V}_{t,T}\rightarrow \mathcal{U}_{t,T}$ such that
for any stopping time $\tau :\Omega \rightarrow \left[ t,T\right] $ and any $%
v_{1},$ $v_{2}\in \mathcal{V}_{t,T}$ with $v_{1}\equiv v_{2}$ on $[[t,\tau
]] $, it holds that $\alpha \left( v_{1}\right) \equiv \alpha \left(
v_{2}\right) $ on $[[t,\tau ]]$. Nonanticipative strategies for player II on
$\left[ t,T\right] $, denoted by $\beta :\mathcal{U}_{t,T}\rightarrow
\mathcal{V}_{t,T}$, are defined similarly. The set of all nonanticipative
strategies $\alpha $ (resp., $\beta $) for player I (resp., II) on $\left[
t,T\right] $ is denoted by $\mathcal{A}_{t,T}$ (resp., $\mathcal{B}_{t,T}$).
(Recall that $[[s,\tau ]]=\left\{ \left( r,\omega \right) \in \left[ 0,T%
\right] \times \Omega ,s\leq r\leq \tau \left( \omega \right) \right\} $).
\end{definition}

Assume that (A1)-(A3) are in force, for any $u\left( \cdot \right) \times
v\left( \cdot \right) \in \mathcal{U}_{t,T}\times \mathcal{V}_{t,T}$, it is
easy to check that FBSDEs (\ref{SDE1})-(\ref{BSDE1}) admit a unique $\mathcal{F}_{t}$%
-adapted strong solution denoted by the triple
\begin{equation*}
(X^{t,x;u,v},Y^{t,x;u,v},Z^{t,x;u,v})\in \mathcal{N}^{2}\left[ 0,T\right]
\end{equation*}
(See Pardoux and Peng \cite{PP1}).

As in Peng in \cite{P3}, given any impulse controls $u\left( \cdot \right)
\times v\left( \cdot \right) \in \mathcal{U}_{t,T}\times \mathcal{V}_{t,T}$,
we introduce the following cost functional:
\begin{equation}
J(t,x;u\left( \cdot \right) ,v\left( \cdot \right) )=\left.
Y_{s}^{t,x;u,v}\right\vert _{s=t},\qquad \left( t,x\right) \in \left[ 0,T%
\right] \times \mathbb{R}^{n}.  \label{c1}
\end{equation}%
Under assumptions (A1)-(A3), the gain functional $J(t,x;u\left( \cdot
\right) ,v\left( \cdot \right) )$, defined by (\ref{c1}) is well defined for
every $(t,x)\in \left[ t,T\right] \times \mathbb{R}^{n}$, $u\in \mathcal{U}%
_{t,T}$, and $v\left( \cdot \right) \in \mathcal{V}_{t,T}$. We are
interested in two \emph{value functions} of the stochastic differential
games of the following type:%
\begin{equation}
V^{-}\left( t,x\right) =\inf_{\beta \in \mathcal{B}_{t,T}}\sup_{u\in
\mathcal{U}_{t,T}}J(t,x;u,\beta \left( u\right) ),\text{ }\left( t,x\right)
\in \left[ 0,T\right] \times \mathbb{R}^{n}  \label{value1}
\end{equation}%
and
\begin{equation}
V^{+}\left( t,x\right) =\sup_{\alpha \in \mathcal{A}_{t,T}}\inf_{v\in
\mathcal{V}_{t,T}}J(t,x;\alpha \left( v\right) ,v),\text{ }\left( t,x\right)
\in \left[ 0,T\right] \times \mathbb{R}^{n}  \label{value2}
\end{equation}%
for every $\left( t,x\right) \in \left[ 0,T\right] \times \mathbb{R}^{n}.$
When $V^{-}=V^{+},$ we say that the game admits a value and $V:=V^{-}=V^{+}$
is called the value function of the game. Since the value function (\ref%
{value1}) and (\ref{value2}) are defined by the solution of controlled BSDE (%
\ref{BSDE1}), $V^{-}\left( V^{+}\right) $ is \emph{well-defined}. Moreover,
they are both bounded $\mathcal{F}_{t}$-measurable random variables.
Nonetheless, we shall prove that $V^{-}\left( V^{+}\right) $ are even
deterministic.

Note that $\inf $ and $\sup $ in this paper should be interpreted via the
essential infimum and the essential supremum with respect to indexed
families of random variables (see Karatzas and Shreve \cite{KSE}). For
reader's convenience, we recall the notion of $essinf$ of processes. Given a
family of real-valued random variables $\eta _{\alpha }$, $\alpha \in I$, a
random variable $\eta $ is said to be $essinf_{\alpha \in I}\eta _{\alpha }$%
, if (i) $\eta \leq \eta _{\alpha }$, $P$-a.s., for any $\alpha \in I$; (ii)
if there is another random variable $\xi $ such that $\xi \leq \eta _{\alpha
}$, $P$-a.s., for any $\alpha \in I$, then $\xi \leq \eta $, $P$-a.s. The
random variable $esssup_{\alpha \in I}\eta _{\alpha }$ can be introduced now
by the relation $esssup_{\alpha \in I}\eta _{\alpha }=-essinf_{\alpha \in
I}\left( -\eta _{\alpha }\right) $. Finally, recall that $essinf_{\alpha \in
I}\eta _{\alpha }=\inf_{n\geq 1}\eta _{\alpha _{n}}$ for some countable
family $\left( \alpha _{n}\right) \subset I$; $esssup_{\alpha \in I}\eta
_{\alpha }$ has the same property. \newline
We need the following estimations for BSDE, whose proof can be seen in
Proposition 3.2 of Briand et al. \cite{BDHPS}.

\begin{lemma}
\label{l1}Let $\left( y^{i},z^{i}\right) ,$ $i=1,2,$ be the solution to the
following
\begin{equation}
y^{i}\left( t\right) =\xi ^{i}+\int_{t}^{T}f^{i}\left( s,y^{i}\left(
s\right) ,z^{i}\left( s\right) \right) \mathrm{d}s-\int_{t}^{T}z^{i}\left(
s\right) \mathrm{d}W_{s},  \label{estbdsde}
\end{equation}%
where $\xi ^{i}\in L^{2}\left( \Omega ,\mathcal{F}_{T},P\right) $ with $%
\mathbb{E}\left[ \left\vert \xi ^{i}\right\vert ^{\beta }\right] <\infty ,$ $%
f^{i}\left( s,y^{i},z^{i}\right) $ satisfies the conditions \emph{(A2)}%
\textsl{,} and
\begin{equation*}
\mathbb{E}\left[ \left( \int_{t}^{T}\left\vert f^{i}\left( s,y^{i}\left(
s\right) ,z^{i}\left( s\right) \right) \right\vert \mathrm{d}s\right)
^{\beta }\right] <\infty .
\end{equation*}%
Then, for some $\beta \geq 2,$ there exists a positive constant $C_{\beta }$
such that
\begin{eqnarray*}
&&\mathbb{E}\left[ \sup_{0\leq t\leq T}\left\vert y^{1}\left( t\right)
-y^{2}\left( t\right) \right\vert ^{\beta }+\left( \int_{0}^{T}\left\vert
z^{1}\left( s\right) -z^{2}\left( s\right) \right\vert ^{2}\mathrm{d}%
s\right) ^{\frac{\beta }{2}}\right] \\
&\leq &C_{\beta }\mathbb{E}\Bigg [\left\vert \xi ^{1}-\xi ^{2}\right\vert
^{\beta }+\left( \int_{t}^{T}\left\vert f^{1}\left( s,y^{1}\left( s\right)
,z^{1}\left( s\right) \right) -f^{2}\left( s,y^{2}\left( s\right)
,z^{2}\left( s\right) \right) \right\vert \mathrm{d}s\right) ^{\beta }\Bigg ]%
.
\end{eqnarray*}%
Particularly, whenever putting $\xi ^{2}=0,$ $f^{2}=0,$ one has%
\begin{equation*}
\mathbb{E}\left[ \sup_{0\leq t\leq T}\left\vert y^{1}\left( t\right)
\right\vert ^{\beta }+\left( \int_{0}^{T}\left\vert z^{1}\left( s\right)
\right\vert ^{2}\mathrm{d}s\right) ^{\frac{\beta }{2}}\right] \leq C_{\beta }%
\mathbb{E}\Bigg [\left\vert \xi ^{1}\right\vert ^{\beta }+\left(
\int_{t}^{T}\left\vert f^{1}\left( s,0,0\right) \right\vert \mathrm{d}%
s\right) ^{\beta }\Bigg ].
\end{equation*}
\end{lemma}

We recall the following well-known comparison theorem (see Barles, Buckdahn,
and Pardoux \cite{BBP}, Proposition 2.6) for BSDE.

\begin{lemma}[Comparison theorem]
\label{comp}Let $\left( y^{i},z^{i}\right) ,$ $i=1,2,$ be the solution to
the following
\begin{equation}
y^{i}\left( t\right) =\xi ^{i}+\int_{t}^{T}f^{i}\left(
s,y_{s}^{i},z_{s}^{i}\right) \mathrm{d}s-\int_{t}^{T}z_{s}^{i}\mathrm{d}%
W_{s},  \label{cbsde}
\end{equation}%
where $\mathbb{E}\left[ \left\vert \xi ^{i}\right\vert ^{2}\right] <\infty ,$
$f^{i}\left( s,y^{i},z^{i}\right) $ satisfies the conditions \emph{(A2), }$%
i=1,2$\textsl{.} Under assumption \emph{(A2)}, BSDE \emph{(\ref{cbsde})}
admits a unique adapted solution $\left( y^{i},z^{i}\right) ,$ respectively,
for $i=1,2$. Furthermore, if (i) $\xi ^{1}\geq \xi ^{2},$ a.s.; (ii) $%
f^{1}\left( t,y,z\right) \geq f^{2}\left( t,y,z\right) ,$ a.e., for any $%
\left( t,y,z\right) \in \left[ 0,T\right] \times \mathbb{R}\times \mathbb{R}%
^{d}.$ Then we have $y^{1}\left( t\right) \geq y^{2}\left( t\right) ,$ a.s.
\end{lemma}

\section{Dynamic Programming Principle}

\label{sec3}

In this section, we present the DPP for our stochastic differential games in
the framework of BSDE. The following lemma announces that the values
functions are deterministic, which is important to investigate the other
properties of value functions.

\begin{lemma}
\label{deter}Let $\left( t,x\right) \in \left[ 0,T\right] \times \mathbb{R}%
^{n}.$ Under assumptions \emph{(A1)-(A3),} $V^{-}\left( t,x\right) =\mathbb{E%
}\left[ V^{-}\left( t,x\right) \right] $, $P$-a.s. Since $V^{-}\left(
t,x\right) $ coincides with its deterministic version $\mathbb{E}\left[
V^{-}\left( t,x\right) \right] ,$ we can consider $V^{-}:\left[ 0,T\right]
\times \mathbb{R}^{n}\rightarrow \mathbb{R}$ as a deterministic function. An
analogous statement holds for the value function $V^{+}.$
\end{lemma}

The proof is displayed in Appendix \ref{APP}.

We shall consider the value functions obtained by no impulse controls, which
is useful to prove the H\"{o}lder continuity of value functions in the
sequel.

\begin{lemma}
\label{noim}Assume that \emph{(A1)-(A3)} are in force, then the lower and
upper value functions are given by%
\begin{equation}
V^{-}\left( t,x\right) =\inf_{\beta \in \mathcal{\bar{B}}_{t,T}}\sup_{u\in
\mathcal{\bar{U}}_{t,T}}J(t,x;u,\beta \left( u\right) ),\text{ }\left(
t,x\right) \in \left[ 0,T\right] \times \mathbb{R}^{n}  \label{vnoim1}
\end{equation}%
and
\begin{equation}
V^{+}\left( t,x\right) =\sup_{\alpha \in \mathcal{\bar{A}}_{t,T}}\inf_{v\in
\mathcal{\bar{V}}_{t,T}}J(t,x;\alpha \left( v\right) ,v),\text{ }\left(
t,x\right) \in \left[ 0,T\right] \times \mathbb{R}^{n},  \label{vnoim2}
\end{equation}%
where $\mathcal{\bar{U}}_{t,T}$ and $\mathcal{\bar{V}}_{t,T}$ contain all
the impulse controls in $\mathcal{U}_{t,T}$ and $\mathcal{V}_{t,T}$,
respectively, which have no impulses at time $t$. Similarly, $\mathcal{\bar{A%
}}_{t,T}$ and $\mathcal{\bar{B}}_{t,T}$ are subsets of $\mathcal{A}_{t,T}$
and $\mathcal{B}_{t,T}$, respectively. In particular, they contain all the
nonanticipative strategies with values in $\mathcal{\bar{U}}_{t,T}$ and $%
\mathcal{\bar{V}}_{t,T}$, respectively.
\end{lemma}

\paragraph{Proof}

We borrow the idea from Cosso's work \cite{co}. Fix $\epsilon >0.$ Let $u\in
\mathcal{U}_{t,T}\backslash \mathcal{\bar{U}}_{t,T}$ and $\beta \in \mathcal{%
B}_{t,T}\backslash \mathcal{\bar{B}}_{t,T}.$ Then, let $v:=\beta \left(
u\right) \in \mathcal{V}_{t,T}$ and $\bar{\beta}\left( \check{u}\right) =%
\bar{v}\in \mathcal{\bar{V}}_{t,T}$ for any $\check{u}\in \mathcal{U}_{t,T}$
for some $\bar{\beta}\in \mathcal{\bar{B}}_{t,T}$. Hence, we have to prove
that there exist $\bar{u}\in \mathcal{\bar{U}}_{t,T}$ and $\bar{v}\in
\mathcal{\bar{V}}_{t,T}$ such that%
\begin{equation*}
\left\vert J\left( t,x;u,v\right) -J\left( t,x;\bar{u},\bar{v}\right)
\right\vert ^{2}\leq \epsilon .
\end{equation*}%
We may suppose $v:=\mathcal{V}_{t,T}\backslash \mathcal{\bar{V}}_{t,T}$; in
the other case can be proved similarly.

At the beginning, let $u$ and $v$ have only a single impulse at time $t$.
Therefore, there exist two $\left[ t,T\right] $-valued $\mathcal{F}$%
-stopping times $\tau $ and $\rho $, with $P\left( \tau =t\right) >0$ and $%
P\left( \rho =t\right) >0$, such that $u=\xi \mathbf{1}_{\left[ \tau ,T%
\right] }+\hat{u}$ and $v=\eta \mathbf{1}_{\left[ \rho ,T\right] }+\hat{v},$
where $\hat{u}=\sum_{m\geq 1}\xi _{m}\mathbf{1}_{\left[ \tau _{m},T\right]
}\in \mathcal{\bar{U}}_{t,T},$ $\hat{v}=\sum_{l\geq 1}\eta _{l}\mathbf{1}_{%
\left[ \rho _{l},T\right] }\in \mathcal{\bar{V}}_{t,T},$ $\xi $ is an $%
\mathcal{F}_{\tau }$-measurable $U$-valued random variable and $\eta $ is an
$\mathcal{F}_{\rho }$-measurable $V$-valued random variable. Let us
introduce the following stopping times:%
\begin{equation*}
\tau _{n}=\left( \tau +\frac{1}{n}\right) \mathbf{1}_{\left\{ \tau
=t\right\} }+\tau \mathbf{1}_{\left\{ \tau >t\right\} }\text{ and }\rho
_{n}=\left( \rho +\frac{1}{n}\right) \mathbf{1}_{\left\{ \rho =t\right\}
}+\rho \mathbf{1}_{\left\{ \rho >t\right\} }
\end{equation*}%
Clearly, $\tau _{n}\rightarrow \tau $ and $\rho _{n}\rightarrow \rho $, as $%
n $ approaches to infinity, $P$-a.s. Now we define the admissible impulse
controls as follows:
\begin{equation*}
u_{n}=\xi \mathbf{1}_{\left[ \tau _{n},T\right] }+\hat{u}\in \mathcal{\bar{U}%
}_{t,T}\text{ and }v_{n}=\eta \mathbf{1}_{\left[ \rho _{n},T\right] }+\hat{v}%
\in \mathcal{\bar{V}}_{t,T}.
\end{equation*}%
By Proposition 3.2 in \cite{BDHPS} and Lemma \ref{deter}, we have the
following estimate:%
\begin{eqnarray*}
&&\left\vert J\left( t,x;u,v\right) -J\left( t,x;\bar{u},\bar{v}\right)
\right\vert ^{2} \\
&=&\left\vert Y_{t}^{t,x;u,v}-Y_{t}^{t,x;u_{n},v_{n}}\right\vert ^{2} \\
&\leq &C\mathbb{E}\Bigg [\Big |\Phi \left( X_{T}^{t,x;u,v}\right) -\Phi
\left( X_{T}^{t,x;u_{n},v_{n}}\right) \\
&&+\sum_{l\geq 1}\chi \left( \rho _{l},\eta _{l}\right) \mathbf{1}_{\left\{
\rho _{l}\leq T\right\} }-\sum_{m\geq 1}c\left( \tau _{m},\xi _{m}\right)
\mathbf{1}_{\left\{ \tau _{m}\leq T\right\} }\prod_{l\geq 1}\mathbf{1}%
_{\left\{ \tau _{m}\neq \rho _{l}\right\} } \\
&&-\sum_{l\geq 1}\chi \left( \rho _{l},\eta _{l}\right) \mathbf{1}_{\left\{
\rho _{l}\leq T\right\} }-\sum_{m\geq 1}c\left( \tau _{m},\xi _{m}\right)
\mathbf{1}_{\left\{ \tau _{m}\leq T\right\} }\prod_{l\geq 1}\mathbf{1}%
_{\left\{ \tau _{m}\neq \rho _{l}\right\} }\Big |^{2} \\
&&+\left( \int_{t}^{T}\left\vert f\left(
s,X_{s}^{t,x;u,v},Y_{s}^{t,x;u,v},Z_{s}^{t,x;u,v}\right) -f\left(
s,X_{s}^{t,x;u_{n},v_{n}},Y_{s}^{t,x;u_{n},v_{n}},Z_{s}^{t,x;u_{n},v_{n}}%
\right) \right\vert \mathrm{d}s\right) ^{2}\Bigg ].
\end{eqnarray*}%
Note that, for every $s\in \left[ t,T\right] $ by by Burkholder-Davis-Gundy
(B-D-G for short, see \cite{RY}) inequality, $%
X_{s}^{t,x;u_{n},v_{n}}\rightarrow X_{s}^{t,x;u,v}$ as $n\rightarrow \infty $%
, $P$-a.s. Therefore, from Gr\"{o}nwall's inequality and the dominated
convergence theorem, there is an integer $N\geq 1$ such that $\left\vert
J\left( t,x;u,v\right) -J\left( t,x;u_{n},v_{n}\right) \right\vert ^{2}\leq
\epsilon .$ As for multiple impulses at time $t$, one can show the same
result by using the assumptions (A3), which actually leads to the case of
the previous one with only a single impulse at time $t$. We thus complete
the proof.\hfill $\Box $

Now we prove the two values functions are bounded.

\begin{proposition}
Assume that assumptions \emph{(A1)-(A3) }are in force, Then, the lower and
upper value functions are bounded.
\end{proposition}

\paragraph{Proof}

We only consider the lower value function; the other case is analogous. Let $%
\varepsilon >0$; then, by the definition of lower value function (\ref%
{value1}), we have, for any $\left( t,x\right) \in \left[ 0,T\right] \times
\mathbb{R}^{n},$ there exists some $\beta _{\varepsilon }\left( u_{0}\right)
=\sum_{l\geq 1}\eta _{l}^{\varepsilon }\left( u_{0}\right) 1_{\left[ \rho
_{l}^{\varepsilon },T\right] }\in \mathcal{V}_{t,T},$ where $u_{0}\in
\mathcal{U}_{t,T}$ denotes the control with no impulses,
\begin{eqnarray*}
V^{-}\left( t,x\right) &=&\inf_{\beta \in \mathcal{B}_{t,T}}\sup_{u\in
\mathcal{U}_{t,T}}J(t,x;u,\beta \left( u\right) )=\inf_{\beta \in \mathcal{B}%
_{t,T}}\sup_{u\in \mathcal{U}_{t,T}}Y_{t}^{t,x;u,\beta \left( u\right) } \\
&=&\inf_{\beta \in \mathcal{B}_{t,T}}\sup_{u\in \mathcal{U}_{t,T}}\mathbb{E}%
\Bigg [\Phi \left( X_{T}^{t,x;u,\beta \left( u\right) }\right)
+\int_{s}^{T}f\left( r,X_{r}^{t,x;u,\beta \left( u\right)
},Y_{r}^{t,x;u,\beta \left( u\right) },Z_{r}^{t,x;u,\beta \left( u\right)
}\right) \mathrm{d}r \\
&&-\int_{s}^{T}Z_{r}^{t,x;u,\beta \left( u\right) }\mathrm{d}%
W_{r}-\sum_{m\geq 1}c\left( \tau _{m},\xi _{m}\right) \mathbf{1}_{\left\{
\tau _{m}\leq T\right\} }\prod_{l\geq 1}\mathbf{1}_{\left\{ \tau _{m}\neq
\rho _{l}\left( u\right) \right\} } \\
&&+\sum_{l\geq 1}\chi \left( \rho _{l}\left( u\right) ,\eta _{l}\left(
u\right) \right) \mathbf{1}_{\left\{ \rho _{l}\left( u\right) \leq T\right\}
}\Bigg ]
\end{eqnarray*}%
\begin{eqnarray*}
&\geq &Y_{t}^{t,x;u_{0},\beta _{\varepsilon }\left( u_{0}\right)
}-\varepsilon \\
&=&\mathbb{E}\Bigg [\Phi \left( X_{T}^{t,x;u_{0},\beta _{\varepsilon }\left(
u_{0}\right) }\right) +\int_{s}^{T}f\left( r,X_{r}^{t,x;u_{0},\beta
_{\varepsilon }\left( u_{0}\right) },Y_{r}^{t,x;u_{0},\beta _{\varepsilon
}\left( u_{0}\right) },Z_{r}^{t,x;u_{0},\beta _{\varepsilon }\left(
u_{0}\right) }\right) \mathrm{d}r \\
&&-\int_{s}^{T}Z_{r}^{t,x;u_{0},\beta _{\varepsilon }\left( u_{0}\right) }%
\mathrm{d}W_{r}+\sum_{l\geq 1}\chi \left( \rho _{l}\left( u\right) ,\eta
_{l}\left( u\right) \right) \mathbf{1}_{\left\{ \rho _{l}\left( u\right)
\leq T\right\} }\Bigg ]-\varepsilon \\
&\geq &\mathbb{E}\Bigg [\Phi \left( X_{T}^{t,x;u_{0},\beta _{\varepsilon
}\left( u_{0}\right) }\right) \\
&&+\int_{s}^{T}f\left( r,X_{r}^{t,x;u_{0},\beta _{\varepsilon }\left(
u_{0}\right) },Y_{r}^{t,x;u_{0},\beta _{\varepsilon }\left( u_{0}\right)
},Z_{r}^{t,x;u_{0},\beta _{\varepsilon }\left( u_{0}\right) }\right) \mathrm{%
d}r\Bigg ]-\varepsilon .
\end{eqnarray*}

The last inequality above is based on the condition (\ref{a1}) and the
comparison theorem (see Proposition 2.6 in \cite{BBP}). Due to $f$ and $\Phi
$ are bounded, we deduce that $V^{-}$ is bounded from below. In a similar
way, we can prove that $V^{-}$ is also bounded from above.\hfill $\Box $

After getting the first result on value functions, we now focus on (the
generalized) DPP in the framework our stochastic differential game (\ref%
{BSDE1}), (\ref{value1}) and (\ref{value2}). To this end, we should first
define the family of (backward) semigroups associated with FBSDEs (\ref%
{BSDE1}). As a matter of fact, this concept of stochastic \emph{backward
semigroups} was first introduced by Peng \cite{P3} which was employed to
investigate the DPP for stochastic control problems.

For every the initial data $\left( t,x\right) \in \left[ 0,T\right] \times
\mathbb{R}^{n}$, a positive number $\delta \leq T-t$, two admissible impulse
control processes $u\in \mathcal{U}_{t,T}$ and $v\in \mathcal{V}_{t,T}$, and
a real-valued random variable $\eta \in L^{2}\left( \Omega ,\mathcal{F}%
_{t+\delta },P;\mathbb{R}\right) $), we define
\begin{equation*}
G_{s,t+\delta }^{t,x;u,v}\left[ \eta +\Theta _{t+\delta }^{u,v}\right] :=%
\mathcal{Y}_{s}^{t,x;u,v},
\end{equation*}%
for $s\in \left[ t,t+\delta \right] $ where
\begin{equation*}
\Theta _{s}^{u,v}:=\sum_{l\geq 1}\chi \left( \rho _{l},\eta _{l}\right)
\mathbf{1}_{\left\{ \rho _{l}\leq s\right\} }-\sum_{m\geq 1}c\left( \tau
_{m},\xi _{m}\right) \mathbf{1}_{\left\{ \tau _{m}\leq s\right\}
}\prod_{l\geq 1}\mathbf{1}_{\left\{ \tau _{m}\neq \rho _{l}\right\} }
\end{equation*}%
the couple $\left( \mathcal{Y}_{s}^{t,x;u,v},\mathcal{Z}_{s}^{t,x;u,v}%
\right) _{t\leq s\leq t+\delta }$ is the solution of the following BSDE with
the time horizon $t+\eta $:%
\begin{equation}
\mathcal{Y}_{s}^{t,x;u,v}=\eta +\Theta _{t+\delta }^{u,v}+\int_{s}^{t+\delta
}f\left( r,X_{r}^{t,x;u,v},\mathcal{Y}_{r}^{t,x;u,v},\mathcal{Z}%
_{r}^{t,x;u,v}\right) \mathrm{d}r-\int_{s}^{t+\delta }\mathcal{Z}%
_{r}^{t,x;u,v}\mathrm{d}W_{r},  \label{BG}
\end{equation}%
for $s\in \left[ t,t+\delta \right] $ and $X^{t,x;u,v}$ is the solution to
SDE (\ref{SDE1}). Then, obviously, for the solution $\left(
Y^{t,x;u,v},Z^{t,x;u,v}\right) $ to BSDE (\ref{BSDE1}), we have%
\begin{equation*}
G_{t,T}^{t,x;u,v}\left[ \Phi \left( X_{T}^{t,x;u,v}\right) +\Theta _{T}^{u,v}%
\right] =G_{t,t+\delta }^{t,x;u,v}\left[ Y_{t+\delta }^{t,x;u,v}+\Theta
_{t+\delta }^{u,v}\right] .
\end{equation*}%
Indeed,
\begin{eqnarray*}
J\left( t,x;u,v\right) &=&Y_{t}^{t,x;u,v}=G_{t,T}^{t,x;u,v}\left[ \Phi
\left( X_{T}^{t,x;u,v}\right) +\Theta _{T}^{u,v}\right] \\
&=&G_{t,t+\delta }^{t,x;u,v}\left[ Y_{t+\delta }^{t,x;u,v}+\Theta _{t+\delta
}^{u,v}\right] \\
&=&G_{t,t+\delta }^{t,x;u,v}\left[ J\left( t+\delta ,X_{t+\delta
}^{t,x;u,v};u,v\right) \right] .
\end{eqnarray*}
Now we are ready to derive the the dynamic programming principle (DPP for
short), by virtue of backward semigroups introduced above, in which the
impulse control can be regarded as a terminal condition. This principle is
important tool to character the viscosity solution of corresponding H-J-B
equation (see Section \ref{sec4}).

\begin{theorem}
\label{dpp}Suppose that assumptions \emph{(A1)-(A3)} hold. Then, the value
function $V^{-}$ admits the following DPP: For any $0\leq t<t+\delta \leq T,$
$x\in \mathbb{R}^{n},$%
\begin{equation*}
V^{-}\left( t,x\right) =\inf_{\beta \in \mathcal{B}_{t,T}}\sup_{u\in
\mathcal{U}_{t,T}}G_{t,t+\delta }^{t,x;u,\beta \left( u\right) }\left[
V^{-}\left( t+\delta ,X_{t+\delta }^{t,x;u,\beta \left( u\right) }\right)
+\Theta _{t+\delta }^{u,\beta \left( u\right) }\right] ,P\text{-a.e.}
\end{equation*}%
An analogous statement holds for the value function $V^{+}.$
\end{theorem}

\paragraph{Proof}

We prove only the dynamic programming principle only for $V^{-}$; the other
case is analogous.

Put
\begin{equation*}
V_{\delta }^{-}\left( t,x\right) =\inf_{\beta \in \mathcal{B}_{t,t+\delta
}}\sup_{u\in \mathcal{U}_{t,t+\delta }}G_{t,t+\delta }^{t,x;u,\beta \left(
u\right) }\left[ V^{-}\left( t+\delta ,X_{t+\delta }^{t,x;u,\beta \left(
u\right) }\right) +\Theta _{t+\delta }^{u,\beta \left( u\right) }\right] .
\end{equation*}%
We proceed the proof that $V_{\delta }^{-}\left( t,x\right) $ coincides with
$V^{-}\left( t,x\right) $ into the following steps.

In the first step, we shall prove $V_{\delta }^{-}\left( t,x\right) \geq
V^{-}\left( t,x\right) .$ To this end, we have
\begin{eqnarray*}
V_{\delta }^{-}\left( t,x\right) &=&\inf_{\beta \in \mathcal{B}_{t,t+\delta
}}\sup_{u\in \mathcal{U}_{t,t+\delta }}G_{t,t+\delta }^{t,x;u,\beta \left(
u\right) }\left[ V^{-}\left( t+\delta ,X_{t+\delta }^{t,x;u,\beta \left(
u\right) }\right) +\Theta _{t+\delta }^{u,\beta \left( u\right) }\right] \\
&=&\inf_{\beta \in \mathcal{B}_{t,t+\delta }}\mathcal{I}_{\delta }\left(
t,x,\beta \right) ,
\end{eqnarray*}%
where the notation $\mathcal{I}_{\delta }\left( t,x,\beta \right)
=\sup_{u\in \mathcal{U}_{t,t+\delta }}\mathcal{I}_{\delta }\left(
t,x,u,\beta \left( u\right) \right) $ with
\begin{equation*}
\mathcal{I}_{\delta }\left( t,x,u,v\right) =G_{t,t+\delta }^{t,x;u,v}\left[
V^{-}\left( t+\delta ,X_{t+\delta }^{t,x;u,v}\right) +\Theta _{t+\delta
}^{u,v}\right] ,\text{ }P\text{-a.s.}
\end{equation*}%
and for some sequences $\left\{ \beta _{i}\right\} _{i\geq 1}\subset
\mathcal{B}_{t,t+\delta }$ such that $V_{\delta }^{-}\left( t,x\right)
=\inf_{i\geq 1}\mathcal{I}_{\delta }\left( t,x,\beta _{i}^{1}\right) ,$ $P$%
-a.s. Let $\varepsilon >0$ and set $\bar{\Upsilon}_{i}:=\left\{ \mathcal{I}%
_{\delta }\left( t,x,\beta _{i}^{1}\right) -\varepsilon \leq V_{\delta
}^{-}\left( t,x\right) \right\} \in \mathcal{F}_{t},$ $i\geq 1.$ Construct $%
\Upsilon _{i}:=\bar{\Upsilon}_{i}\backslash \left( \cup _{k=1}^{i-1}\bar{%
\Upsilon}_{k}\right) .$ Certainly, $\left\{ \Upsilon _{i}\right\} _{i\geq 1}$
forms an $\left( \Omega ,\mathcal{F}\right) $-partition, moreover, $\beta
^{\varepsilon ,1}:=\sum_{i\geq 1}\mathbf{1}_{\Upsilon _{i}}\beta _{i}^{1}\in
\mathcal{B}_{t,t+\delta }.$ According the existence and uniquenss of FBSDEs (%
\ref{BSDE1}), it follows that $\mathcal{I}_{\delta }\left( t,x,u^{1},\beta
^{\varepsilon ,1}\left( u^{1}\right) \right) =\sum_{i\geq 1}\mathbf{1}%
_{\Upsilon _{i}}\mathcal{I}_{\delta }\left( t,x,u,\beta _{i}^{1}\left(
u\right) \right) ,$ $P$-a.s., for each $u^{1}\in \mathcal{U}_{t,t+\delta }.$
Next, for $\forall u^{1}\in \mathcal{U}_{t,t+\delta }$
\begin{eqnarray}
V_{\delta }^{-}\left( t,x\right) &\geq &\sum_{i\geq 1}\mathbf{1}_{\Upsilon
_{i}}\mathcal{I}_{\delta }\left( t,x,\beta _{i}^{1}\right) -\varepsilon
\notag \\
&\geq &\sum_{i\geq 1}\mathbf{1}_{\Upsilon _{i}}\mathcal{I}_{\delta }\left(
t,x,u^{1},\beta _{i}^{1}\right) -\varepsilon  \notag \\
&=&\mathcal{I}_{\delta }\left( t,x,u^{1},\beta ^{\varepsilon ,1}\right)
-\varepsilon  \notag \\
&=&G_{t,t+\delta }^{t,x;u^{1},\beta ^{\varepsilon ,1}\left( u^{1}\right) }
\left[ V^{-}\left( t+\delta ,X_{t+\delta }^{t,x;u^{1},\beta ^{\varepsilon
,1}\left( u^{1}\right) }\right) +\Theta _{t+\delta }^{u^{1},\beta
^{\varepsilon ,1}\left( u^{1}\right) }\right] -\varepsilon ,\text{ }P\text{%
-a.s.}.  \label{dpp0}
\end{eqnarray}%
We now focus on the time interval $\left[ t+\delta ,T\right] .$ From the
definition of $V_{\delta }^{-}\left( t,x\right) $, we also deduce that, with
help of previous idea and Lemma \ref{noim}, for any $y\in \mathbb{R}^{n}$,
there exists $\beta _{y}^{\varepsilon }\in \mathcal{\bar{B}}_{t+\delta ,T}$
for each $u^{2}\in \mathcal{U}_{t+\delta ,T}$ such that%
\begin{equation}
V^{-}\left( t+\delta ,y\right) \geq \sup_{u^{2}\in \mathcal{U}_{t+\delta
,T}}J\left( t+\delta ,y,u^{2},\beta _{y}^{\varepsilon }\left( u^{2}\right)
\right) -\varepsilon ,\text{ }P\text{-a.s.}  \label{dpp1}
\end{equation}%
Now consider a decomposition of $\mathbb{R}^{n}$, namely, $\sum_{i\geq 1}%
\mathcal{O}_{i}=\mathbb{R}^{n}$ such that diam$\left( \mathcal{O}_{i}\right)
\leq \varepsilon ,$ for each $i\geq 1.$ Take any $y_{i}\in \mathcal{O}_{i}$
fixed, $i\geq 1$ and define $\overline{X_{t+\delta }^{t,x;u^{1},\beta
^{\varepsilon ,1}\left( u^{1}\right) }}=\sum_{i\geq 1}y_{i}\mathbf{1}%
_{\left\{ X_{t+\delta }^{t,x;u^{1},\beta ^{\varepsilon ,1}\left(
u^{1}\right) }\in \mathcal{O}_{i}\right\} }.$ Clearly, we always have $%
\left\vert X_{t+\delta }^{t,x;u^{1},\beta ^{\varepsilon ,1}\left(
u^{1}\right) }-\overline{X_{t+\delta }^{t,x;u^{1},\beta ^{\varepsilon
,1}\left( u^{1}\right) }}\right\vert \leq \varepsilon ,$ almost on $\Omega $%
, for each $u^{1}\in \mathcal{U}_{t,t+\delta }.$ For every $y_{i}\in
\mathcal{O}_{i},$ one can seek $\beta _{y_{i}}^{\varepsilon }\in \mathcal{%
\bar{B}}_{t+\delta ,T}$ such that (\ref{dpp1}) holds true.

We introduce the strategy $\beta _{u^{1}}^{2,\varepsilon }\in \mathcal{\bar{B%
}}_{t+\delta ,T}$ as follows:
\begin{equation}
\beta _{u^{1}}^{2,\varepsilon }:=\sum_{i\geq 1}\mathbf{1}_{\left\{
X_{t+\delta }^{t,x;u^{1},\beta ^{\varepsilon ,1}\left( u^{1}\right) }\in
\mathcal{O}_{i}\right\} }\beta _{y_{i}}^{\varepsilon }\in \mathcal{\bar{B}}%
_{t+\delta ,T}.  \label{cs1}
\end{equation}%
Set $\beta _{y_{i}}^{\varepsilon }\left( u^{2}\right) =\sum_{l\geq 1}\eta
_{l}^{i}\left( u^{2}\right) \mathbf{1}_{\left[ \rho _{l}^{i}\left(
u^{2}\right) ,T\right] },$ for $u^{2}\in \mathcal{U}_{t+\delta ,T};$ then
\begin{equation*}
\beta _{u^{1}}^{2,\varepsilon }\left( u^{2}\right) :=\sum_{l\geq 1}\eta
_{l,u^{1}}^{2,\varepsilon }\left( u^{2}\right) \mathbf{1}_{\left[ \rho
_{l,u^{1}}^{2,\varepsilon }\left( u^{2}\right) ,T\right] },
\end{equation*}%
where
\begin{equation}
\eta _{l,u^{1}}^{2,\varepsilon }\left( u^{2}\right) :=\sum_{i\geq 1}\mathbf{1%
}_{\left\{ X_{t+\delta }^{t,x;u^{1},\beta ^{\varepsilon ,1}\left(
u^{1}\right) }\in \mathcal{O}_{i}\right\} }\eta _{l}^{i}\left( u^{2}\right)
\label{cs2}
\end{equation}%
and
\begin{equation}
\rho _{l,u^{1}}^{2,\varepsilon }\left( u^{2}\right) :=\sum_{i\geq 1}\mathbf{1%
}_{\left\{ X_{t+\delta }^{t,x;u^{1},\beta ^{\varepsilon ,1}\left(
u^{1}\right) }\in \mathcal{O}_{i}\right\} }\rho _{l}^{i}\left( u^{2}\right) .
\label{cs3}
\end{equation}%
It is possible to define a new strategy $\beta ^{\varepsilon }\left(
u\right) $ from $\beta ^{\varepsilon ,1}\left( u^{1}\right) \in \mathcal{B}%
_{t,t+\delta }$\ and $\beta _{u^{1}}^{2,\varepsilon }\left( u^{2}\right) \in
\mathcal{\bar{B}}_{t+\delta ,T}$ where $u^{1}=u_{\left[ t,t+\delta \right]
}, $ $u^{2}=u_{\left( t+\delta ,T\right] }$ (see Definition \ref{s2})$,$ in
the following way:

Let
\begin{equation}
\beta ^{\varepsilon ,1}\left( u^{1}\right) =\sum_{l\geq 1}\eta
_{l}^{1,\varepsilon }\left( u^{1}\right) \mathbf{1}_{\left[ \rho
_{l}^{1,\varepsilon }\left( u^{1}\right) ,T\right] },\text{ }\beta
_{u^{1}}^{\varepsilon ,2}\left( u^{2}\right) =\sum_{l\geq 1}\eta
_{l,u^{1}}^{2,\varepsilon }\left( u^{2}\right) \mathbf{1}_{\left[ \rho
_{l,u^{1}}^{2,\varepsilon }\left( u^{2}\right) ,T\right] ,}  \label{cs4}
\end{equation}%
then $\beta ^{\varepsilon }\left( u\right) =\sum_{l\geq 1}\eta
_{l}^{\varepsilon }\left( u\right) \mathbf{1}_{\left[ \rho _{l}^{\varepsilon
}\left( u\right) ,T\right] },$ where
\begin{equation}
\eta _{l}^{\varepsilon }\left( u\right) =\eta _{l}^{1,\varepsilon }\left(
u^{1}\right) \mathbf{1}_{\left\{ l\leq \mu _{t,t+\delta }\left( \beta
^{\varepsilon ,1}\left( u^{1}\right) \right) \right\} }+\eta _{l-\mu
_{t,t+\delta }\left( \beta ^{\varepsilon ,1}\left( u\right) \right)
,u^{1}}^{2,\varepsilon }\left( u^{2}\right) \mathbf{1}_{\left\{ l>\mu
_{t,t+\delta }\left( \beta ^{\varepsilon ,1}\left( u^{1}\right) \right)
\right\} }  \label{cs5}
\end{equation}%
and%
\begin{equation}
\rho _{l}^{\varepsilon }\left( u\right) =\rho _{l}^{1,\varepsilon }\left(
u^{1}\right) \mathbf{1}_{\left\{ l\leq \mu _{t,t+\delta }\left( \beta
^{\varepsilon ,1}\left( u^{1}\right) \right) \right\} }+\rho _{l-\mu
_{t,t+\delta }\left( \beta ^{\varepsilon ,1}\left( u\right) \right)
,u^{1}}^{2,\varepsilon }\left( u^{2}\right) \mathbf{1}_{\left\{ l>\mu
_{t,t+\delta }\left( \beta ^{\varepsilon ,1}\left( u\right) \right) \right\}
},  \label{cs6}
\end{equation}%
where $\mu _{t,t+\delta }$ is defined in (\ref{fimp}).

Next we shall show that $\beta ^{\varepsilon }\left( u\right) $ is
nonanticipating: Indeed, let $\kappa :\Omega \rightarrow \left[ t,T\right] $
be an $\mathcal{F}$-stopping time and $u,$ $u^{\prime }\in \mathcal{U}_{t,T}$
be such that $u\equiv u^{\prime }$ on $\left[ t,\kappa \right] $.
Decomposing $u,$ $u^{\prime }$ into $u_{1},$ $u_{1}^{\prime }\in \mathcal{B}%
_{t,t+\delta },$ $u_{2},$ $u_{2}^{\prime }\in \mathcal{\bar{B}}_{t+\delta
,T} $ such that $u=u_{1}\oplus u_{2},$ $u^{\prime }=u_{1}^{\prime }\oplus
u_{2}^{\prime }$ where $u_{1}\oplus u_{2}$ (the same for $u_{1}^{\prime
}\oplus u_{2}^{\prime }$) is defined as follows:

Let
\begin{equation}
u_{1}=\sum_{l\geq 1}\eta _{l}^{1}\mathbf{1}_{\left[ \rho _{l}^{1},T\right]
},u_{2}=\sum_{l\geq 1}\eta _{l}^{2}\mathbf{1}_{\left[ \rho _{l}^{2},T\right]
}.  \label{lt1}
\end{equation}%
then $u_{1}\oplus u_{2}=\sum_{l\geq 1}\eta _{l}^{\oplus }\mathbf{1}_{\left[
\rho _{l}^{\oplus },T\right] },$ where
\begin{equation}
\eta _{l}^{\oplus }=\eta _{l}^{1}\mathbf{1}_{\left\{ l\leq \mu _{t,t+\delta
}\left( u_{1}\right) \right\} }+\eta _{l-\mu _{t,t+\delta }\left(
u_{1}\right) }^{2}\mathbf{1}_{\left\{ l>\mu _{t,t+\delta }\left(
u_{1}\right) \right\} }  \label{lt2}
\end{equation}%
and%
\begin{equation}
\rho _{l}^{\oplus }=\rho _{l}^{1}\mathbf{1}_{\left\{ l\leq \mu _{t,t+\delta
}\left( u_{1}\right) \right\} }+\rho _{l-\mu _{t,t+\delta }\left(
u_{1}\right) }^{2}\mathbf{1}_{\left\{ l>\mu _{t,t+\delta }\left(
u_{1}\right) \right\} }.  \label{lt3}
\end{equation}%
Thus $\beta ^{\varepsilon }\left( u\right) =\beta ^{\varepsilon ,1}\left(
u^{1}\right) \oplus \beta _{u^{1}}^{2,\varepsilon }\left( u^{2}\right) $ for
$u^{1}=u_{\left[ t,t+\delta \right] },$ $u^{2}=u_{\left( t,T\right] }.$ We
immediately have $\beta ^{\varepsilon ,1}\left( u_{1}\right) =\beta
^{\varepsilon ,1}\left( u_{1}^{\prime }\right) $ since $u_{1}=u_{1}^{\prime
} $ on $\left[ t,\kappa \wedge t+\delta \right] $. On the other hand, $%
u_{2}=u_{2}^{\prime }$ on $\left( t+\delta ,\kappa \vee t+\delta \right] $
and on $\left\{ \kappa >t+\delta \right\} $, we have $X_{t+\delta
}^{t,x;u_{1},\beta ^{\varepsilon ,1}\left( u_{1}\right) }=X_{t+\delta
}^{t,x;u_{1}^{\prime },\beta ^{\varepsilon ,1}\left( u_{1}^{\prime }\right)
}.$ This yields our desired result.

Fix $u\in \mathcal{U}_{t,T}$ arbitrarily and decompose into $u^{1}=u_{\left[
t,t+\delta \right] },$ $u^{2}=u_{\left( t+\delta ,T\right] }$. Then, from (%
\ref{dpp0}) and Proposition 2.6 in \cite{BBP}, we obtain%
\begin{eqnarray}
&&V_{\delta }^{-}\left( t,x\right)  \notag \\
&\geq &G_{t,t+\delta }^{t,x;u^{1},\beta ^{\varepsilon ,1}\left( u^{1}\right)
}\left[ V^{-}\left( t+\delta ,X_{t+\delta }^{t,x;u^{1},\beta ^{\varepsilon
,1}\left( u^{1}\right) }\right) +\Theta _{t+\delta }^{u^{1},\beta
^{\varepsilon ,1}\left( u^{1}\right) }\right] -\varepsilon  \notag \\
&\geq &G_{t,t+\delta }^{t,x;u^{1},\beta ^{\varepsilon ,1}\left( u^{1}\right)
}\left[ V^{-}\left( t+\delta ,\overline{X_{t+\delta }^{t,x;u^{1},\beta
^{\varepsilon ,1}\left( u^{1}\right) }}\right) +\Theta _{t+\delta
}^{u^{1},\beta ^{\varepsilon ,1}\left( u^{1}\right) }\right] -C\varepsilon
\notag \\
&=&G_{t,t+\delta }^{t,x;u^{1},\beta ^{\varepsilon ,1}\left( u^{1}\right) }
\left[ \sum_{i\geq 1}\mathbf{1}_{\left\{ X_{t+\delta }^{t,x;u^{1},\beta
^{\varepsilon ,1}\left( u^{1}\right) }\in \mathcal{O}_{i}\right\}
}V^{-}\left( t+\delta ,y_{i}\right) +\Theta _{t+\delta }^{u^{1},\beta
^{\varepsilon ,1}\left( u^{1}\right) }\right] -C\varepsilon ,\text{ }P\text{%
-a.s.}  \label{b1}
\end{eqnarray}%
From (\ref{b1}) and Proposition 2.6 in \cite{BBP}, it follows%
\begin{eqnarray*}
&&V_{\delta }^{-}\left( t,x\right) \\
&\geq &G_{t,t+\delta }^{t,x;u^{1},\beta ^{\varepsilon ,1}\left( u^{1}\right)
}\left[ \sum_{i\geq 1}\mathbf{1}_{\left\{ X_{t+\delta }^{t,x;u^{1},\beta
^{\varepsilon ,1}\left( u^{1}\right) }\in \mathcal{O}_{i}\right\} }J\left(
t+\delta ,y_{i},u^{2},\beta _{y_{i}}^{\varepsilon }\left( u^{2}\right)
\right) +\Theta _{t+\delta }^{u^{1},\beta ^{\varepsilon ,1}\left(
u^{1}\right) }\right] -C\varepsilon \\
&=&G_{t,t+\delta }^{t,x;u^{1},\beta ^{\varepsilon ,1}\left( u^{1}\right) }
\left[ J\left( t+\delta ,\overline{X_{t+\delta }^{t,x;u^{1},\beta
^{\varepsilon ,1}\left( u^{1}\right) }},u^{2},\beta _{y_{i}}^{\varepsilon
}\left( u^{2}\right) \right) +\Theta _{t+\delta }^{u^{1},\beta ^{\varepsilon
,1}\left( u^{1}\right) }\right] -C\varepsilon \\
&\geq &G_{t,t+\delta }^{t,x;u^{1},\beta ^{\varepsilon ,1}\left( u^{1}\right)
}\left[ J\left( t+\delta ,X_{t+\delta }^{t,x;u^{1},\beta ^{\varepsilon
,1}\left( u^{1}\right) },u^{2},\beta _{y_{i}}^{\varepsilon }\left(
u^{2}\right) \right) +\Theta _{t+\delta }^{u^{1},\beta ^{\varepsilon
,1}\left( u^{1}\right) }\right] -C\varepsilon \\
&=&G_{t,t+\delta }^{t,x;u,\beta ^{\varepsilon }\left( u\right) }\left[
Y_{t+\delta }^{t,x;u,\beta ^{\varepsilon }\left( u\right) }\right]
-C\varepsilon \\
&=&Y_{t}^{t,x;u,\beta ^{\varepsilon }\left( u\right) }-C\varepsilon ,\text{ }%
P\text{-a.s., for every }u\in \mathcal{U}_{t,T}.
\end{eqnarray*}%
Therefore, we obtain%
\begin{eqnarray*}
V_{\delta }^{-}\left( t,x\right) &\geq &\sup_{u\in \mathcal{U}_{t,T}}J\left(
t,x;u;\beta ^{\varepsilon }\left( u\right) \right) -C\varepsilon \\
&\geq &\inf_{\beta \in \mathcal{B}_{t,T}}\sup_{u\in \mathcal{U}%
_{t,T}}J\left( t,x;u;\beta ^{\varepsilon }\left( u\right) \right)
-C\varepsilon \\
&=&V^{-}\left( t,x\right) -C\varepsilon ,\text{ }P\text{-a.s.}
\end{eqnarray*}%
We now deal with the other case: $V_{\delta }^{-}\left( t,x\right) \leq
V^{-}\left( t,x\right) .$

Let $\beta \in \mathcal{B}_{t,T}$ be arbitrarily chosen and $u_{2}\in
\mathcal{\bar{U}}_{t+\delta ,T}$. Define the restriction of $\beta $ to $%
\mathcal{U}_{t,t+\delta }$ as $\beta ^{1}\left( u_{1}\right) :=\beta \left(
u_{1}\oplus u_{2}\right) _{\left[ t,t+\delta \right] },$ $u_{1}\in \mathcal{U%
}_{t,t+\delta }.$ The nonanticipativity property of $\beta $ indicates that $%
\beta ^{1}$ is independent of the special choice of $u_{2}\in \mathcal{\bar{U%
}}_{t+\delta ,T}$. From the definition of $V_{\delta }^{-}\left( t,x\right)
, $
\begin{equation*}
V_{\delta }^{-}\left( t,x\right) \leq \sup_{u_{1}\in \mathcal{U}_{t,t+\delta
}}G_{t,t+\delta }^{t,x;u_{1},\beta ^{1}\left( u_{1}\right) }\left[
V^{-}\left( t+\delta ,X_{t+\delta }^{t,x;u_{1},\beta ^{1}\left( u_{1}\right)
}\right) \right] ,\text{ }P\text{-a.s.}
\end{equation*}%
Consider $\mathcal{I}_{\delta }\left( t,x,\beta ^{1}\right) =\sup_{u^{1}\in
\mathcal{U}_{t,t+\delta }}\mathcal{I}_{\delta }\left( t,x,u^{1},\beta
^{1}\left( u^{1}\right) \right) ;$ then there exists a sequence $\left\{
u_{i}^{1}\right\} _{i\geq 1}\subset \mathcal{\bar{U}}_{t,t+\delta }$ such
that $\mathcal{I}_{\delta }\left( t,x,\beta ^{1}\right) =\sup_{i\geq 1}%
\mathcal{I}_{\delta }\left( t,x,u_{i}^{1},\beta ^{1}\left( u_{i}^{1}\right)
\right) ,$ $P$-a.s. With the same technique as before, for any $\varepsilon
>0,$ set $\bar{\Lambda}_{i}:=\left\{ \mathcal{I}_{\delta }\left( t,x,\beta
_{i}^{1}\right) \leq \mathcal{I}_{\delta }\left( t,x,u_{i}^{1},\beta
^{1}\left( u_{i}^{1}\right) \right) +\varepsilon \right\} \in \mathcal{F}%
_{t},$ $i\geq 1.$ Construct $\Lambda _{i}:=\bar{\Lambda}_{i}\backslash
\left( \cup _{k=1}^{i-1}\bar{\Lambda}_{k}\right) .$ Certainly, $\left\{
\Lambda _{i}\right\} _{i\geq 1}$ forms an $\left( \Omega ,\mathcal{F}\right)
$-partition, moreover, $u_{1}^{\varepsilon }:=\sum_{i\geq 1}1_{\Lambda
_{i}}u_{i}^{1}\in \mathcal{U}_{t,t+\delta }.$ From the existence and
uniqueness of FBSDEs (\ref{BSDE1}), we deduce that $\mathcal{I}_{\delta
}\left( t,x,u_{1}^{\varepsilon },\beta ^{1}\left( u_{1}^{\varepsilon
}\right) \right) =\sum_{i\geq 1}\mathcal{I}_{\delta }\left(
t,x,u_{i}^{1},\beta ^{1}\left( u_{i}^{1}\right) \right) ,$ $P$-a.s. Then,
\begin{eqnarray*}
V_{\delta }^{-}\left( t,x\right) &\leq &\mathcal{I}_{\delta }\left(
t,x,\beta ^{1}\right) \leq \sum_{i\geq 1}\mathbf{1}_{\Lambda _{i}}\mathcal{I}%
_{\delta }\left( t,x,u_{i}^{1},\beta ^{1}\left( u_{i}^{1}\right) \right)
+\varepsilon \\
&=&\mathcal{I}_{\delta }\left( t,x,u_{1}^{\varepsilon },\beta ^{1}\left(
u_{1}^{\varepsilon }\right) \right) +\varepsilon \\
&=&G_{t,t+\delta }^{t,x;u_{1}^{\varepsilon },\beta ^{1}\left(
u_{1}^{\varepsilon }\right) }\left[ V^{-}\left( t+\delta ,X_{t+\delta
}^{t,x;u_{1}^{\varepsilon },\beta ^{1}\left( u_{1}^{\varepsilon }\right)
}\right) +\Theta _{t+\delta }^{u_{1}^{\varepsilon },\beta ^{1}\left(
u_{1}^{\varepsilon }\right) }\right] +\varepsilon ,\text{ }P\text{-a.s.}
\end{eqnarray*}%
Noting that $\beta ^{1}\left( \cdot \right) :=\beta \left( \cdot \oplus
u_{2}\right) \in \mathcal{B}_{t,t+\delta }$ does not depend on $u_{2}\in
\mathcal{\bar{U}}_{t+\delta ,T}$, we can construct $\beta ^{2}\left(
u_{2}\right) :=\beta \left( u_{1}^{\varepsilon }\oplus u_{2}\right) _{\left[
t+\delta ,T\right] },$ for each $u_{2}\in \mathcal{\bar{U}}_{t+\delta ,T}$
such that $\beta ^{2}:\mathcal{\bar{U}}_{t+\delta ,T}\rightarrow \mathcal{%
\bar{V}}_{t+\delta ,T}$ belongs to $\mathcal{\bar{B}}_{t+\delta ,T}$, due to
$\beta \in \mathcal{B}_{t,T}$. Therefore, from the definition of $%
V^{-}\left( t+\delta ,y\right) $ and Lemma \ref{noim}$,$ we have, for any $%
y\in \mathbb{R}^{n}$,
\begin{equation*}
V^{-}\left( t+\delta ,X_{t+\delta }^{t,x;u_{1}^{\varepsilon },\beta
^{1}\left( u_{1}^{\varepsilon }\right) }\right) \leq \sup_{u_{2}\in \mathcal{%
U}_{t+\delta ,T}}J\left( t+\delta ,X_{t+\delta }^{t,x;u_{1}^{\varepsilon
},\beta ^{1}\left( u_{1}^{\varepsilon }\right) };u_{2},\beta ^{2}\left(
u_{2}\right) \right) .
\end{equation*}%
There exists a sequence $\left\{ u_{2}^{i}\right\} _{i\geq 1}\subset
\mathcal{U}_{t+\delta ,T}$ such that
\begin{eqnarray*}
&&\sup_{u_{2}\in \mathcal{U}_{t+\delta ,T}}J\left( t+\delta ,X_{t+\delta
}^{t,x;u_{1}^{\varepsilon },\beta ^{1}\left( u_{1}^{\varepsilon }\right)
};u_{2},\beta ^{2}\left( u_{2}\right) \right) \\
&=&\sup_{i\geq 1}J\left( t+\delta ,X_{t+\delta }^{t,x;u_{1}^{\varepsilon
},\beta ^{1}\left( u_{1}^{\varepsilon }\right) };u_{2}^{i},\beta ^{2}\left(
u_{2}^{i}\right) \right) .
\end{eqnarray*}%
Then with the same technique as before, for any $\varepsilon >0,$ set
\begin{eqnarray*}
\bar{\Pi}_{i} &:&=\Big \{\sup_{u_{2}\in \mathcal{U}_{t+\delta ,T}}J\left(
t+\delta ,X_{t+\delta }^{t,x;u_{1}^{\varepsilon },\beta ^{1}\left(
u_{1}^{\varepsilon }\right) };u_{2},\beta ^{2}\left( u_{2}\right) \right) \\
&\leq &J\left( t+\delta ,X_{t+\delta }^{t,x;u_{1}^{\varepsilon },\beta
^{1}\left( u_{1}^{\varepsilon }\right) };u_{2}^{i},\beta ^{2}\left(
u_{2}^{i}\right) \right) +\varepsilon \Big\}\in \mathcal{F}_{t+\delta
},i\geq 1.
\end{eqnarray*}
Construct $\Pi _{i}:=\bar{\Pi}_{i}\backslash \left( \cup _{k=1}^{i-1}\bar{\Pi%
}_{k}\right) .$ Certainly, $\left\{ \Pi _{i}\right\} _{i\geq 1}$ also forms
an $\left( \Omega ,\mathcal{F}\right) $-partition, moreover, $%
u_{2}^{\varepsilon }:=\sum_{i\geq 1}\mathbf{1}_{\Pi _{i}}u_{2}^{i}\in
\mathcal{U}_{t+\delta ,T}.$ Then, $\beta ^{2}\left( u_{2}^{\varepsilon
}\right) =\sum_{j\geq 1}\mathbf{1}_{\Pi _{i}}.$ We construct a new strategy $%
\beta \left( u_{1}^{\varepsilon }\oplus u_{2}^{\varepsilon }\right) =\beta
^{1}\left( u_{1}^{\varepsilon }\right) \oplus \beta ^{2}\left(
u_{2}^{\varepsilon }\right) .$ From the existence and uniqueness of FBSDEs (%
\ref{BSDE1}), we have%
\begin{eqnarray}
J\left( t+\delta ,X_{t+\delta }^{t,x;u_{1}^{\varepsilon },\beta ^{1}\left(
u_{1}^{\varepsilon }\right) };u_{2}^{\varepsilon },\beta ^{2}\left(
u_{2}^{\varepsilon }\right) \right) &=&Y_{t+\delta }^{t+\delta ,X_{t+\delta
}^{t,x;u_{1}^{\varepsilon },\beta ^{1}\left( u_{1}^{\varepsilon }\right)
};u_{2}^{\varepsilon },\beta ^{2}\left( u_{2}^{\varepsilon }\right) }  \notag
\\
&=&\sum_{j\geq 1}\mathbf{1}_{\Pi _{i}}Y_{t+\delta }^{t+\delta ,X_{t+\delta
}^{t,x;u_{1}^{\varepsilon },\beta ^{1}\left( u_{1}^{\varepsilon }\right)
};u_{2}^{j},\beta ^{2}\left( u_{2}^{j}\right) }  \notag \\
&=&\sum_{j\geq 1}\mathbf{1}_{\Pi _{i}}J\left( t+\delta ,X_{t+\delta
}^{t,x;u_{1}^{\varepsilon },\beta ^{1}\left( u_{1}^{\varepsilon }\right)
};u_{2}^{j},\beta ^{2}\left( u_{2}^{j}\right) \right) .  \label{dd1}
\end{eqnarray}%
Therefore,
\begin{eqnarray}
V^{-}\left( t+\delta ,X_{t+\delta }^{t,x;u_{1}^{\varepsilon },\beta
^{1}\left( u_{1}^{\varepsilon }\right) }\right) &\leq &\sup_{u_{2}\in
\mathcal{U}_{t+\delta ,T}}J\left( t+\delta ,X_{t+\delta
}^{t,x;u_{1}^{\varepsilon },\beta ^{1}\left( u_{1}^{\varepsilon }\right)
};u_{2},\beta ^{2}\left( u_{2}\right) \right)  \notag \\
&\leq &\sum_{j\geq 1}\mathbf{1}_{\Pi _{i}}Y_{t+\delta
}^{t,x;u_{1}^{\varepsilon }\oplus u_{2}^{j},\beta \left( u_{1}^{\varepsilon
}\oplus u_{2}^{j}\right) }+\varepsilon  \notag \\
&=&Y_{t+\delta }^{t,x;u_{1}^{\varepsilon }\oplus u_{2}^{\varepsilon },\beta
\left( u_{1}^{\varepsilon }\oplus u_{2}^{\varepsilon }\right) }+\varepsilon
\notag \\
&=&Y_{t+\delta }^{t,x;u^{\varepsilon },\beta \left( u^{\varepsilon }\right)
}+\varepsilon ,  \label{dd2}
\end{eqnarray}%
where $u^{\varepsilon }=u_{1}^{\varepsilon }\oplus u_{2}^{\varepsilon }\in
\mathcal{U}_{t,T}.$ Repeating the method before, from (\ref{dd1}) and (\ref%
{dd2}) and Proposition 2.6 in \cite{BBP}, we have%
\begin{eqnarray*}
V_{\delta }^{-}\left( t,x\right) &\leq &G_{t,t+\delta
}^{t,x;u_{1}^{\varepsilon },\beta ^{1}\left( u_{1}^{\varepsilon }\right) }
\left[ Y_{t+\delta }^{t,x;u^{\varepsilon },\beta \left( u^{\varepsilon
}\right) }+\Theta _{t+\delta }^{u^{\varepsilon },\beta \left( u^{\varepsilon
}\right) }\right] +C\varepsilon \\
&=&G_{t,t+\delta }^{t,x;u^{\varepsilon },\beta ^{1}\left( u^{\varepsilon
}\right) }\left[ Y_{t+\delta }^{t,x;u^{\varepsilon },\beta \left(
u^{\varepsilon }\right) }+\Theta _{t+\delta }^{u^{\varepsilon },\beta \left(
u^{\varepsilon }\right) }\right] +C\varepsilon \\
&=&Y_{t}^{t,x;u^{\varepsilon },\beta \left( u^{\varepsilon }\right)
}+C\varepsilon \\
&\leq &\sup_{u\in \mathcal{U}_{t,T}}Y_{t}^{t,x;u,\beta \left( u\right)
}+C\varepsilon ,\text{ }P\text{-a.s.,}
\end{eqnarray*}%
which holds for all $\beta \in \mathcal{B}_{t,T}.$
\begin{equation*}
V_{\delta }^{-}\left( t,x\right) \leq \inf_{\beta \in \mathcal{B}%
_{t,T}}\sup_{u\in \mathcal{U}_{t,T}}Y_{t}^{t,x;u,\beta \left( u\right)
}+C\varepsilon =V^{-}\left( t,x\right) +C\varepsilon
\end{equation*}%
Now letting $\varepsilon \rightarrow 0,$ we get the desired result, $%
V_{\delta }^{-}\left( t,x\right) \leq V^{-}\left( t,x\right) $. The proof is
completed.\hfill $\Box $

Next, we will show that the continuity of value functions with respect to $x$
and $t$. Due to the influence of Brownian motion, the value function will
proved to be Lipschitz continuity on $x$, but $1/2$ h\"{o}lder on $t$.

\begin{proposition}
\label{vlip}Assume that assumptions \emph{(A1)-(A3) }are in force. Then the
lower value function $V^{-}\left( t,x\right) $ is $\frac{1}{2}$-H\"{o}lder
continuous in $t$: There exists a constant $C>0$ such that, for every $%
\left( t,x\right) \in \left[ 0,T\right) \times \mathbb{R}^{n}$
\begin{eqnarray}
\left\vert V^{-}\left( t,x\right) -V^{-}\left( t,x^{\prime }\right)
\right\vert &\leq &C\left\vert x-x^{\prime }\right\vert ,  \label{lipx} \\
\left\vert V^{-}\left( t,x\right) -V^{-}\left( t^{\prime },x\right)
\right\vert &\leq &C\left\vert t-t^{\prime }\right\vert ^{\frac{1}{2}}.
\label{lipt}
\end{eqnarray}
\end{proposition}

\paragraph{Proof}

The first property of the lower value function $V^{-}\left( t,x\right) $
which we present is an immediate consequence of Proposition 3.2 in \cite%
{BDHPS}.

Let $\left( t,x\right) \in \left[ 0,T\right) \times \mathbb{R}^{n}$ and $%
\delta >0$ be arbitrarily given such that $0<\delta <T-t$. Then for every $%
\varepsilon >0$, thanks to Lemma \ref{noim}, there exist $u_{\varepsilon
}\in \mathcal{U}_{t,T}$ and $\beta _{\varepsilon }\in \mathcal{\bar{B}}%
_{t,T} $ (This ensures no impulse on initial state) such that%
\begin{eqnarray}
&&V^{-}\left( t,x\right) -V^{-}\left( t+\delta ,x\right)  \notag \\
&\leq &G_{t,T}^{t,x;u_{\varepsilon },\hat{\beta}_{\varepsilon }\left(
u_{\varepsilon }\right) }\left( \Phi \left( X_{T}^{t,x;u_{\varepsilon },\hat{%
\beta}_{\varepsilon }\left( u_{\varepsilon }\right) }\right) \right)
-G_{t+\delta ,T}^{t,x;\hat{u}_{\varepsilon },\beta _{\varepsilon }\left(
\hat{u}_{\varepsilon }\right) }\left( \Phi \left( X_{T}^{t,x;\hat{u}%
_{\varepsilon },\beta _{\varepsilon }\left( \hat{u}_{\varepsilon }\right)
}\right) \right) +\varepsilon  \label{holder1}
\end{eqnarray}%
where $\hat{u}_{\varepsilon }\in \mathcal{U}_{t+\delta ,T}$ and $\hat{\beta}%
_{\varepsilon }\in \mathcal{B}_{t,T}$ will be determined soon. Indeed, from (%
\ref{vnoim1}) and (\ref{vnoim2}), there exist $u_{\varepsilon }\in \mathcal{U%
}_{t,T}$ and $\beta _{\varepsilon }\in \mathcal{\bar{B}}_{t,T}$ such that

\begin{equation}
G_{t,T}^{t,x;u_{\varepsilon },\hat{\beta}_{\varepsilon }\left(
u_{\varepsilon }\right) }\left( \Phi \left( X_{T}^{t,x;u_{\varepsilon },\hat{%
\beta}_{\varepsilon }\left( u_{\varepsilon }\right) }\right) \right) \geq
V^{-}\left( t,x\right) -\frac{\varepsilon }{2}  \label{t1}
\end{equation}%
and%
\begin{equation}
G_{t+\delta ,T}^{t,x;\hat{u}_{\varepsilon },\beta _{\varepsilon }\left( \hat{%
u}_{\varepsilon }\right) }\left( \Phi \left( X_{T}^{t,x;\hat{u}_{\varepsilon
},\beta _{\varepsilon }\left( \hat{u}_{\varepsilon }\right) }\right) \right)
\leq V^{-}\left( t+\delta ,x\right) +\frac{\varepsilon }{2}.  \label{t2}
\end{equation}%
From (\ref{t1}) and (\ref{t2}), one can obtain (\ref{holder1}) easily.

We postulate that $u_{\varepsilon }=\sum_{m\geq 1}\xi _{m}^{\varepsilon }%
\mathbf{1}_{\left[ \tau _{m}^{\varepsilon },T\right] }\in \mathcal{U}_{t,T}$%
; now define $\hat{u}_{\varepsilon }$ as $\hat{u}_{\varepsilon }=\sum_{\tau
_{m}\leq t+\delta }\xi _{m}^{\varepsilon }\mathbf{1}_{t+\delta }+\sum_{\tau
_{m}>t+\delta }\xi _{m}^{\varepsilon }\mathbf{1}_{\left[ \tau
_{m}^{\varepsilon },T\right] }.$ Observe that $\hat{u}_{\varepsilon }$ is an
impulse control constructing from $u_{\varepsilon }$ via gathering all the
impulses in the interval $\left[ t,t+\delta \right] .$ To eliminate the
impulses on time $t+\delta $ of player II$,$ we define $v_{\varepsilon }=%
\hat{\beta}_{\varepsilon }\left( u\right) =\beta _{\varepsilon }\left( \hat{u%
}_{\varepsilon }\right) \in \mathcal{\bar{V}}_{t+\delta ,T}$ for any $u\in
\mathcal{U}_{t,T}.$ After this work, (\ref{holder1}) can be written as
\begin{eqnarray}
&&V^{-}\left( t,x\right) -V^{-}\left( t+\delta ,x\right)  \notag \\
&\leq &G_{t,T}^{t,x;u_{\varepsilon },v_{\varepsilon }}\left( \Phi \left(
X_{T}^{t,x;u_{\varepsilon },v_{\varepsilon }}\right) \right) -G_{t+\delta
,T}^{t,x;\hat{u}_{\varepsilon },v_{\varepsilon }}\left( \Phi \left(
X_{T}^{t,x;\hat{u}_{\varepsilon },v_{\varepsilon }}\right) \right)
+\varepsilon .  \label{holder2}
\end{eqnarray}%
We deal with%
\begin{eqnarray}
&&G_{t,T}^{t,x;u_{\varepsilon },v_{\varepsilon }}\left( \Phi \left(
X_{T}^{t,x;u_{\varepsilon },v_{\varepsilon }}\right) \right) -G_{t+\delta
,T}^{t,x;\hat{u}_{\varepsilon },v_{\varepsilon }}\left( \Phi \left(
X_{T}^{t,x;\hat{u}_{\varepsilon },v_{\varepsilon }}\right) \right)  \notag \\
&=&\Phi \left( T,X_{T}^{t,x;u_{\varepsilon },v_{\varepsilon }}\right) -\Phi
\left( X_{T}^{t,x;\hat{u}_{\varepsilon },v_{\varepsilon }}\right) \mathrm{d}%
r+\int_{s}^{T}Z_{r}^{t,x;\hat{u}_{\varepsilon },v_{\varepsilon }}\mathrm{d}%
W_{r}-\int_{s}^{T}Z_{r}^{t,x;u_{\varepsilon },v_{\varepsilon }}\mathrm{d}%
W_{r}  \notag \\
&&+\left[ c\left( t+\delta ,\sum_{\tau _{m}\leq t+\delta }\xi
_{m}^{\varepsilon }\right) \mathbf{1}_{t+\delta }+\sum_{\tau _{m}>t+\delta
}c\left( \tau _{m}^{\varepsilon },\xi _{m}^{\varepsilon }\right) \mathbf{1}_{%
\left[ \tau _{m}^{\varepsilon },T\right] }\right] \mathbf{1}_{\left\{ \tau
_{m}^{\varepsilon }\leq T\right\} }\prod_{l\geq 1}\mathbf{1}_{\left\{ \tau
_{m}^{\varepsilon }\neq \rho _{l}^{\varepsilon }\right\} }  \notag \\
&&-\sum_{m\geq 1}c\left( \tau _{m}^{\varepsilon },\xi _{m}^{\varepsilon
}\right) \mathbf{1}_{\left\{ \tau _{m}^{\varepsilon }\leq T\right\}
}\prod_{l\geq 1}\mathbf{1}_{\left\{ \tau _{m}^{\varepsilon }\neq \rho
_{l}^{\varepsilon }\right\} }  \notag \\
&&+\int_{s}^{T}f\left( r,X_{r}^{t,x;u_{\varepsilon },v_{\varepsilon
}},Y_{r}^{t,x;u_{\varepsilon },v_{\varepsilon }},Z_{r}^{t,x;u_{\varepsilon
},v_{\varepsilon }}\right)  \notag \\
&&-\int_{s}^{T}f\left( r,X_{r}^{t,x;\hat{u}_{\varepsilon },v_{\varepsilon
}},Y_{r}^{t,x;\hat{u}_{\varepsilon },v_{\varepsilon }},Z_{r}^{t,x;\hat{u}%
_{\varepsilon },v_{\varepsilon }}\right) \mathrm{d}r,  \label{holder3}
\end{eqnarray}%
but from conditions (\ref{a2}) and (\ref{a4}), we have%
\begin{equation*}
c\left( t+\delta ,\sum_{\tau _{m}\leq t+\delta }\xi _{m}^{\varepsilon
}\right) \leq \sum_{\tau _{m}\leq t+\delta }c\left( \tau _{m}^{\varepsilon
},\xi _{m}^{\varepsilon }\right) \mathbf{1}_{t+\delta }.
\end{equation*}%
Hence, (\ref{holder3}) yields
\begin{eqnarray}
&&G_{t,T}^{t,x;u_{\varepsilon },v_{\varepsilon }}\left( \Phi \left(
X_{T}^{t,x;u_{\varepsilon },v_{\varepsilon }}\right) \right) -G_{t+\delta
,T}^{t,x;\hat{u}_{\varepsilon },v_{\varepsilon }}\left( \Phi \left(
X_{T}^{t,x;\hat{u}_{\varepsilon },v_{\varepsilon }}\right) \right)  \notag \\
&\leq &\Phi \left( T,X_{T}^{t,x;u_{\varepsilon },v_{\varepsilon }}\right)
-\Phi \left( X_{T}^{t,x;\hat{u}_{\varepsilon },v_{\varepsilon }}\right)
\mathrm{d}r+\int_{s}^{T}Z_{r}^{t,x;\hat{u}_{\varepsilon },v_{\varepsilon }}%
\mathrm{d}W_{r}-\int_{s}^{T}Z_{r}^{t,x;u_{\varepsilon },v_{\varepsilon }}%
\mathrm{d}W_{r}  \notag \\
&&+\int_{s}^{T}f\left( r,X_{r}^{t,x;u_{\varepsilon },v_{\varepsilon
}},Y_{r}^{t,x;u_{\varepsilon },v_{\varepsilon }},Z_{r}^{t,x;u_{\varepsilon
},v_{\varepsilon }}\right)  \notag \\
&&-\int_{s}^{T}f\left( r,X_{r}^{t,x;\hat{u}_{\varepsilon },v_{\varepsilon
}},Y_{r}^{t,x;\hat{u}_{\varepsilon },v_{\varepsilon }},Z_{r}^{t,x;\hat{u}%
_{\varepsilon },v_{\varepsilon }}\right) \mathrm{d}r.  \label{holder4}
\end{eqnarray}%
Taking the expectation on both sides of (\ref{holder4}) and noting Lemma \ref%
{deter}, we have%
\begin{eqnarray}
&&G_{t,T}^{t,x;u_{\varepsilon },v_{\varepsilon }}\left( \Phi \left(
X_{T}^{t,x;u_{\varepsilon },v_{\varepsilon }}\right) \right) -G_{t+\delta
,T}^{t,x;\hat{u}_{\varepsilon },v_{\varepsilon }}\left( \Phi \left(
X_{T}^{t,x;\hat{u}_{\varepsilon },v_{\varepsilon }}\right) \right)  \notag \\
&\leq &\mathbb{E}\Bigg [\Phi \left( X_{T}^{t,x;u_{\varepsilon
},v_{\varepsilon }}\right) -\Phi \left( X_{T}^{t,x;\hat{u}_{\varepsilon
},v_{\varepsilon }}\right) \mathrm{d}r  \notag \\
&&+\int_{s}^{T}f\left( r,X_{r}^{t,x;u_{\varepsilon },v_{\varepsilon
}},Y_{r}^{t,x;u_{\varepsilon },v_{\varepsilon }},Z_{r}^{t,x;u_{\varepsilon
},v_{\varepsilon }}\right)  \notag \\
&&-\int_{s}^{T}f\left( r,X_{r}^{t,x;\hat{u}_{\varepsilon },v_{\varepsilon
}},Y_{r}^{t,x;\hat{u}_{\varepsilon },v_{\varepsilon }},Z_{r}^{t,x;\hat{u}%
_{\varepsilon },v_{\varepsilon }}\right) \mathrm{d}r\Bigg ].  \label{holder5}
\end{eqnarray}%
Set $\hat{\Xi}_{r}=\Xi _{r}^{t,x;u_{\varepsilon },v_{\varepsilon }}-\Xi
_{r}^{t,x;\hat{u}_{\varepsilon },v_{\varepsilon }}$, $r\in \left[ t+\delta ,T%
\right] $ for $\Xi =X,$ $Y,$ $Z$. By B-D-G inequality and classical method,
we have%
\begin{equation*}
\mathbb{E}\left[ \sup_{t+\delta \leq r\leq T}\left\vert \hat{X}%
_{r}\right\vert ^{2}+\sup_{t+\delta \leq r\leq T}\left\vert \hat{Y}%
_{r}\right\vert ^{2}+\int_{s}^{T}\left\vert \hat{Z}_{r}\right\vert ^{2}%
\mathrm{d}r\right] \leq C\left\vert t-t^{\prime }\right\vert ^{\frac{1}{2}}.
\end{equation*}%
Therefore,
\begin{eqnarray*}
&&V^{-}\left( t,x\right) -V^{-}\left( t+\delta ,x\right) \\
&\leq &G_{t,T}^{t,x;u_{\varepsilon },v_{\varepsilon }}\left( \Phi \left(
X_{T}^{t,x;u_{\varepsilon },v_{\varepsilon }}\right) \right) -G_{t+\delta
,T}^{t,x;\hat{u}_{\varepsilon },v_{\varepsilon }}\left( \Phi \left(
X_{T}^{t,x;\hat{u}_{\varepsilon },v_{\varepsilon }}\right) \right) \\
&\leq &C\left\vert t-t^{\prime }\right\vert ^{\frac{1}{2}}+\varepsilon .
\end{eqnarray*}%
Letting $\varepsilon \rightarrow 0,$ we get the desired result. We thus
complete the proof.\hfill $\Box $

Now we are concerned on a special case of DPP, that is $s=t$, thanks to
conditions (A3), the multiple impulses can be neglected. It will be useful
in proving that the two value functions are viscosity solutions to the
associated HJBI equation and deriving the so called lower and upper
obstacles. Whilst, it announces that our games problems can interpreted via
optimal stopping times.

\begin{lemma}
\label{st}Assume assumptions \emph{(A1)-(A3)} are in force. Given any $%
\left( t,x\right) \in \left[ 0,T\right] \times \mathbb{R}^{n},$ we have%
\begin{eqnarray}
V^{-}\left( t,x\right) &=&\inf_{\rho \in \mathcal{T}_{t,+\infty }\eta \in
\mathcal{F}_{\rho }}\sup_{\tau \in \mathcal{T}_{t,+\infty },\xi \in \mathcal{%
F}_{\tau }}G_{t,t}^{t,x;u,\beta \left( u\right) }\Big [-c\left( t,\xi
\right) \mathbf{1}_{\left\{ \tau =t\right\} }\mathbf{1}_{\left\{ \rho
=+\infty \right\} }  \notag \\
&&+\chi \left( t,\eta \right) \mathbf{1}_{\left\{ \rho =t\right\}
}+V^{-}\left( t,X_{t}^{t,x;\xi \mathbf{1}_{\left[ \tau ,T\right] },\eta
\mathbf{1}_{\left[ \rho ,T\right] }}\right) \Big ],  \label{st0}
\end{eqnarray}%
where $\mathcal{T}_{t,+\infty }$ is the set of $\mathcal{F}$-stopping times
with values in $\left\{ t,+\infty \right\} $, $\tau \in \mathcal{T}%
_{t,+\infty },$ $\xi \in \mathcal{F}_{\tau },$ $u=\xi \mathbf{1}_{\left[
\tau ,T\right] }$ and $\rho \in \mathcal{T}_{t,+\infty },$ $\eta \in
\mathcal{F}_{\rho },$ $\beta \left( u\right) =\eta \mathbf{1}_{\left[ \rho ,T%
\right] }.$ An analogous statement holds for the upper value function $%
V^{+}\left( t,x\right) $.
\end{lemma}

\paragraph{Proof}

In Theorem \ref{dpp}, consider $V^{-}$ with $\delta =0$:%
\begin{eqnarray*}
V^{-}\left( t,x\right) &=&\inf_{\beta \in \mathcal{B}_{t,T}}\sup_{u\in
\mathcal{U}_{t,T}}G_{t,t}^{t,x;u,\beta \left( u\right) }\Bigg [V^{-}\left(
t,X_{t}^{t,x;u,\beta \left( u\right) }\right) +\sum_{l\geq 1}\chi \left(
\rho _{l},\eta _{l}\right) \mathbf{1}_{\left\{ \rho _{l}=t\right\} } \\
&&-\sum_{m\geq 1}c\left( \tau _{m},\xi _{m}\right) \mathbf{1}_{\left\{ \tau
_{m}=t\right\} }\prod_{l\geq 1}\mathbf{1}_{\left\{ \tau _{m}\neq \rho
_{l}\right\} }\Bigg ].
\end{eqnarray*}%
Given any $u\in \mathcal{U}_{t,T},$ consider the strategy $\beta \left(
u\right) =\eta \mathbf{1}_{\left[ \rho ,T\right] }.$ Let $u=\sum_{m\geq
1}\xi _{m}\mathbf{1}_{\left[ \tau _{m},T\right] },$ then constrcut a new
control $\bar{u}=\xi \mathbf{1}_{\left[ \tau ,T\right] },$ where
\begin{equation*}
\tau =t\left( 1-\prod_{m\geq 1}\mathbf{1}_{\left\{ \tau _{m}>t\right\}
}\right) +\infty \prod_{m\geq 1}\mathbf{1}_{\left\{ \tau _{m}>t\right\} },%
\text{ }\xi =\sum_{m\geq 1}\xi _{m}\mathbf{1}_{\left\{ \tau _{m}=t\right\} }.
\end{equation*}%
Apparently, $\tau \in \mathcal{T}_{t,+\infty }$ and $\xi \in \mathcal{F}%
_{\tau }.$ Meanwhile, we deduce that $X_{t}^{t,x;u,\beta \left( u\right)
}=X_{t}^{t,x;\xi \mathbf{1}_{\left[ \tau ,T\right] },\eta \mathbf{1}_{\left[
\rho ,T\right] }},$ $P$-a.s. By means of (A3), it follows that
\begin{eqnarray*}
&&G_{t,t}^{t,x;u,\beta \left( u\right) }\Big [-\sum_{m\geq 1}c\left( t,\xi
_{m}\right) \mathbf{1}_{\left\{ \tau _{m}=t\right\} }\mathbf{1}_{\left\{
\rho =+\infty \right\} }+\chi \left( t,\eta \right) \mathbf{1}_{\left\{ \rho
=t\right\} }+V^{-}\left( t,X_{t}^{t,x;u,\beta \left( u\right) }\right) \Big ]
\\
&\leq &G_{t,t}^{t,x;u,\beta \left( u\right) }\Big [-c\left( t,\xi \right)
\mathbf{1}_{\left\{ \tau =t\right\} }\mathbf{1}_{\left\{ \rho =+\infty
\right\} }+\chi \left( t,\eta \right) \mathbf{1}_{\left\{ \rho =t\right\}
}+V^{-}\left( t,X_{t}^{t,x;\xi \mathbf{1}_{\left[ \tau ,T\right] },\eta
\mathbf{1}_{\left[ \rho ,T\right] }}\right) \Big ].
\end{eqnarray*}%
As a result, we have
\begin{eqnarray*}
V^{-}\left( t,x\right) &\leq &\inf_{\rho \in \mathcal{T}_{t,+\infty }\eta
\in \mathcal{F}_{\rho }}\sup_{\tau \in \mathcal{T}_{t,+\infty },\xi \in
\mathcal{F}_{\tau }}\mathbb{E}\Big [-c\left( t,\xi \right) \mathbf{1}%
_{\left\{ \tau =t\right\} }\mathbf{1}_{\left\{ \rho =+\infty \right\} } \\
&&+\chi \left( t,\eta \right) \mathbf{1}_{\left\{ \rho =t\right\}
}+V^{-}\left( t,X_{t}^{t,x;\xi \mathbf{1}_{\left[ \tau ,T\right] },\eta
\mathbf{1}_{\left[ \rho ,T\right] }}\right) \Big ].
\end{eqnarray*}%
The reverse inequality can be proved in the analogous way. We end the
proof.\hfill $\Box $

In order to prove the the two value functions satisfy, in the viscosity
sense, the terminal condition to the HJBI equation. We need a useful
technical lemma.

\begin{lemma}
\label{ter}Assume assumption \emph{(A1)-(A3)} are in force. Given any $%
\left( t,x\right) \in \left[ 0,T\right] \times \mathbb{R}^{n},$ we have%
\begin{eqnarray}
&&V^{-}\left( t,x\right)  \notag \\
&=&\inf_{\rho \in \mathcal{T}_{t,+\infty }\eta \in \mathcal{F}_{\rho
}}\sup_{\tau \in \mathcal{T}_{t,+\infty },\xi \in \mathcal{F}_{\tau
}}G_{t,t}^{t,x;u,\beta \left( u\right) }\Big [\Big (-c\left( t,\xi \right)
\mathbf{1}_{\left\{ \tau =t\right\} }\mathbf{1}_{\left\{ \rho =+\infty
\right\} }  \notag \\
&&+\chi \left( t,\eta \right) \mathbf{1}_{\left\{ \rho =t\right\}
}+V^{-}\left( t,X_{t}^{t,x;\xi \mathbf{1}_{\left[ \tau ,T\right] },\eta
\mathbf{1}_{\left[ \rho ,T\right] }}\right) \Big )\left( 1-\mathbf{1}%
_{\left\{ \tau =+\infty ,\rho =+\infty \right\} }\right)  \notag \\
&&+\left( \int_{t}^{T}f\left(
s,X_{s}^{t,x;u_{0},v_{0}},Y_{s}^{t,x;u_{0},v_{0}},Z_{s}^{t,x;u_{0},v_{0}}%
\right) \mathrm{d}s+\Phi \left( X_{T}^{t,x;u_{0},v_{0}}\right) \right)
\mathbf{1}_{\left\{ \tau =+\infty ,\rho =+\infty \right\} }\Big ],
\label{ter1}
\end{eqnarray}%
where $u_{0},$ $v_{0}$ are the controls with no impulses. An analogous
statement holds for the upper value function $V^{+}$.
\end{lemma}

\paragraph{Proof}

For any $\varepsilon >0,$ from the definition of inf, there exist $\rho
^{\varepsilon ,1}\in \mathcal{T}_{t,+\infty },$ $\eta ^{\varepsilon ,1}\in
\mathcal{F}_{\rho ^{\varepsilon ,1}}$ such that
\begin{eqnarray}
&&\text{the right side of}(\ref{ter1})  \notag \\
&\geq &G_{t,t}^{t,x;u,\beta \left( u\right) }\Big [\Big (-c\left( t,\xi
\right) \mathbf{1}_{\left\{ \tau =t\right\} }\mathbf{1}_{\left\{ \rho
^{\varepsilon ,1}=+\infty \right\} }  \notag \\
&&+\chi \left( t,\eta ^{\varepsilon ,1}\right) \mathbf{1}_{\left\{ \rho
^{\varepsilon ,1}=t\right\} }+V^{-}\left( t,X_{t}^{t,x;\xi \mathbf{1}_{\left[
\tau ,T\right] },\eta ^{\varepsilon ,1}\mathbf{1}_{\left[ \rho ^{\varepsilon
,1},T\right] }}\right) \Big )\left( 1-\mathbf{1}_{\left\{ \tau =+\infty
,\rho ^{\varepsilon ,1}=+\infty \right\} }\right)  \notag \\
&&+\left( \int_{t}^{T}f\left(
s,X_{s}^{t,x;u_{0},v_{0}},Y_{s}^{t,x;u_{0},v_{0}},Z_{s}^{t,x;u_{0},v_{0}}%
\right) \mathrm{d}s+\Phi \left( X_{T}^{t,x;u_{0},v_{0}}\right) \right)
\mathbf{1}_{\left\{ \tau =+\infty ,\rho ^{\varepsilon ,1}=+\infty \right\} }%
\Big ]  \notag \\
&&-\varepsilon .  \label{ter2}
\end{eqnarray}%
To deal with $V^{-}\left( t,X_{t}^{t,x;\xi \mathbf{1}_{\left[ \tau ,T\right]
},\eta ^{\varepsilon ,1}\mathbf{1}_{\left[ \rho ^{\varepsilon ,1},T\right]
}}\right) ,$ let $\check{u}\in \mathcal{\bar{U}}_{t,T}.$ From Theorem \ref%
{dpp}, there exsits a strategy $\beta ^{\varepsilon ,2}\in \mathcal{\bar{B}}%
_{t,T}$ such that%
\begin{equation*}
V^{-}\left( t,X_{t}^{t,x;\xi \mathbf{1}_{\left[ \tau ,T\right] },\eta
^{\varepsilon ,1}\mathbf{1}_{\left[ \rho ^{\varepsilon ,1},T\right]
}}\right) \geq \mathbb{E}\left[ J\left( t,X_{t}^{t,x;\xi \mathbf{1}_{\left[
\tau ,T\right] },\eta ^{\varepsilon ,1}\mathbf{1}_{\left[ \rho ^{\varepsilon
,1},T\right] }},\check{u},\beta ^{\varepsilon ,2}\left( \check{u}\right)
\right) \right] -\varepsilon .
\end{equation*}%
\hfill

Define
\begin{eqnarray*}
u^{\varepsilon } &=&\left[ \left( \xi \mathbf{1}_{\left\{ \tau =t\right\}
}+v_{0}\mathbf{1}_{\left\{ \tau =+\infty \right\} }\right) \mathbf{1}_{t}+%
\check{u}\right] \left( 1-\mathbf{1}_{\left\{ \tau =+\infty ,\rho
^{\varepsilon ,1}=+\infty \right\} }\right) +u_{0}\mathbf{1}_{\left\{ \tau
=+\infty ,\rho ^{\varepsilon ,1}=+\infty \right\} }, \\
\beta ^{\varepsilon } &=&\left[ \left( \eta ^{\varepsilon ,1}\mathbf{1}%
_{\left\{ \rho ^{\varepsilon ,1}=t\right\} }+u_{0}\mathbf{1}_{\left\{ \tau
=+\infty \right\} }\right) \mathbf{1}_{t}+\beta ^{\varepsilon ,2}\right]
\left( 1-\mathbf{1}_{\left\{ \tau =+\infty ,\rho ^{\varepsilon ,1}=+\infty
\right\} }\right) +v_{0}\mathbf{1}_{\left\{ \tau =+\infty ,\rho
^{\varepsilon ,1}=+\infty \right\} }.
\end{eqnarray*}%
It is easy to check that $u^{\varepsilon }\in \mathcal{U}_{t,T}$ and $\beta
^{\varepsilon }\in \mathcal{B}_{t,T}.$ Hence, from (\ref{ter1}), (\ref{ter2}%
) and Theorem \ref{dpp}, we have%
\begin{eqnarray*}
&&\inf_{\rho \in \mathcal{T}_{t,+\infty }\eta \in \mathcal{F}_{\rho
}}\sup_{\tau \in \mathcal{T}_{t,+\infty },\xi \in \mathcal{F}_{\tau
}}G_{t,t}^{t,x;u,\beta \left( u\right) }\Big [\Big (-c\left( t,\xi \right)
\mathbf{1}_{\left\{ \tau =t\right\} }\mathbf{1}_{\left\{ \rho =+\infty
\right\} } \\
&&+\chi \left( t,\eta \right) \mathbf{1}_{\left\{ \rho =t\right\}
}+V^{-}\left( t,X_{t}^{t,x;\xi \mathbf{1}_{\left[ \tau ,T\right] },\eta
\mathbf{1}_{\left[ \rho ,T\right] }}\right) \Big )\left( 1-\mathbf{1}%
_{\left\{ \tau =+\infty ,\rho =+\infty \right\} }\right) \\
&&+\left( \int_{t}^{T}f\left(
s,X_{s}^{t,x;u_{0},v_{0}},Y_{s}^{t,x;u_{0},v_{0}},Z_{s}^{t,x;u_{0},v_{0}}%
\right) \mathrm{d}s+\Phi \left( X_{T}^{t,x;u_{0},v_{0}}\right) \right)
\mathbf{1}_{\left\{ \tau =+\infty ,\rho =+\infty \right\} }\Big ] \\
&\geq &G_{t,t}^{t,x;u,\beta \left( u\right) }\Big [\Big (-c\left( t,\xi
\right) \mathbf{1}_{\left\{ \tau =t\right\} }\mathbf{1}_{\left\{ \rho
^{\varepsilon ,1}=+\infty \right\} } \\
&&+\chi \left( t,\eta ^{\varepsilon ,1}\right) \mathbf{1}_{\left\{ \rho
^{\varepsilon ,1}=t\right\} }+J\left( t,X_{t}^{t,x;\xi \mathbf{1}_{\left[
\tau ,T\right] },\eta ^{\varepsilon ,1}\mathbf{1}_{\left[ \rho ^{\varepsilon
,1},T\right] }},\check{u},\beta ^{\varepsilon ,2}\left( \check{u}\right)
\right) \Big )\left( 1-\mathbf{1}_{\left\{ \tau =+\infty ,\rho ^{\varepsilon
,1}=+\infty \right\} }\right) \\
&&+\left( \int_{t}^{T}f\left(
s,X_{s}^{t,x;u_{0},v_{0}},Y_{s}^{t,x;u_{0},v_{0}},Z_{s}^{t,x;u_{0},v_{0}}%
\right) \mathrm{d}s+\Phi \left( X_{T}^{t,x;u_{0},v_{0}}\right) \right)
\mathbf{1}_{\left\{ \tau =+\infty ,\rho ^{\varepsilon ,1}=+\infty \right\} }%
\Big ]-2\varepsilon \\
&=&J\left( t,x;u^{\varepsilon },\beta ^{\varepsilon }\left( u^{\varepsilon
}\right) \right) -2\varepsilon .
\end{eqnarray*}%
Leting $\varepsilon \rightarrow 0,$ we get
\begin{eqnarray*}
V^{-}\left( t,x\right) &\leq &\inf_{\rho \in \mathcal{T}_{t,+\infty }\eta
\in \mathcal{F}_{\rho }}\sup_{\tau \in \mathcal{T}_{t,+\infty },\xi \in
\mathcal{F}_{\tau }}G_{t,t}^{t,x;u,\beta \left( u\right) }\Big [\Big (%
-c\left( t,\xi \right) \mathbf{1}_{\left\{ \tau =t\right\} }\mathbf{1}%
_{\left\{ \rho =+\infty \right\} } \\
&&+\chi \left( t,\eta \right) \mathbf{1}_{\left\{ \rho =t\right\}
}+V^{-}\left( t,X_{t}^{t,x;\xi \mathbf{1}_{\left[ \tau ,T\right] },\eta
\mathbf{1}_{\left[ \rho ,T\right] }}\right) \Big )\left( 1-\mathbf{1}%
_{\left\{ \tau =+\infty ,\rho =+\infty \right\} }\right) \\
&&+\left( \int_{t}^{T}f\left(
s,X_{s}^{t,x;u_{0},v_{0}},Y_{s}^{t,x;u_{0},v_{0}},Z_{s}^{t,x;u_{0},v_{0}}%
\right) \mathrm{d}s+\Phi \left( X_{T}^{t,x;u_{0},v_{0}}\right) \right)
\mathbf{1}_{\left\{ \tau =+\infty ,\rho =+\infty \right\} }\Big ].
\end{eqnarray*}%
The reverse part can be obtained in the same way. We complete the
proof.\hfill $\Box $

\section{HJBI equation: Viscosity approach}

\label{sec4} In the stochastic optimal control theory, the value function is
a solution to the corresponding Hamilton-Jacobi-Bellman equation (H-J-B in
short) whenever it has sufficient regularity (Fleming and Soner \cite{FSon},
Krylov \cite{Krlov}). In other word, it requires that the HJB equation admit
classical solutions, meaning that the solutions be smooth enough (to the
order of derivatives involved in the equation). Unfortunately, this is not
necessarily the case even for some very simple situations. In the stochastic
environment where the diffusion is possibly degenerate, the HJB equation may
in general have no classical solutions either. To overcome this difficulty,
Crandall and Lions introduced the so-called \textit{viscosity solutions} in
the early 1980s (see also \cite{CIL}). This new notion is a kind of
nonsmooth solutions (the value function is continuous, then, the value
function is a solution to the H-J-B equation in the viscosity sense) to
partial differential equations, whose key feature is to replace the
conventional derivatives by the (set-valued) super-/subdifferentials while
maintaining the uniqueness of solutions under very mild conditions. These
make the theory a powerful tool in tackling optimal control problems.

In this section, we consider the following HJBI equation associated to our
stochastic differential games, in which lead to be the same expression for
the two value functions since the two players can not operate at the same
time in the systems, is described by%
\begin{equation}
\left\{
\begin{array}{l}
\max \left\{ V-\mathcal{H}_{\sup }^{\chi }V,\min \left\{ -\frac{\partial }{%
\partial t}V\left( t,x\right) -H\left( t,x,V,DV,D^{2}V\right) ,V-\mathcal{H}%
_{\inf }^{c}V\right\} \right\} =0, \\
V\left( T,x\right) =\Phi \left( x\right) ,\text{ }\left( t,x\right) \in %
\left[ 0,T\right) \times \mathbb{R}^{n},%
\end{array}%
\right.  \label{HJBI}
\end{equation}%
where associated with the Hamiltonians:
\begin{equation}
H\left( t,x,y,p,Q\right) =\left\langle b\left( t,x\right) ,p\right\rangle +%
\frac{1}{2}\text{tr}\left( \sigma \sigma ^{\top }\left( t,x\right) Q\right)
+f\left( t,x,y,p^{\top }\sigma \left( t,x\right) \right)  \label{H}
\end{equation}

and the nonlocal operators $\mathcal{H}_{\sup }^{\chi }V$ and $\mathcal{H}%
_{\inf }^{c}V$ are defined by%
\begin{eqnarray*}
\mathcal{H}_{\sup }^{\chi }V\left( t,x\right) &=&\sup_{y\in U}\left[ V\left(
t,x+y\right) -c\left( t,y\right) \right] , \\
\mathcal{H}_{\inf }^{c}V\left( t,x\right) &=&\inf_{z\in V}\left[ V\left(
t,x+z\right) +\chi \left( t,z\right) \right] ,
\end{eqnarray*}%
for any $\left( t,x\right) \in \left[ 0,T\right) \times \mathbb{R}^{n},$ $%
y\in \mathbb{R},$ $p\in \mathbb{R}^{n},$ $Q\in \mathbb{S}^{n}$ where $%
\mathbb{S}^{n}$ denotes the set of $n\times n$ symmetric matrices. The
coefficients $b,$ $\sigma ,$ $f$, $\Phi $, $\chi $ and $c$ are supposed to
satisfy (A1)-(A3).

We next prove that the lower value function $V\left( t,x\right) $ introduced
by (\ref{HJBI}) is the viscosity solution of (\ref{HJBI}). We extend Cosso's
work \cite{co} for stochastic differential games involving impulse controls
into Peng's BSDE's framework. The difficulties related with this extension
come from the fact that now, contrarily to the framework of stochastic
control theory studied by Peng, we have to do with stochastic differential
games in which strategies are played versus controls. In order to overcome
these difficulties in the proof that $V^{-}$ is a viscosity supersolution,
we have, in particular, to enrich Peng's BSDE method. On the other hand, the
proof that $V^{-}$ is a viscosity subsolution is not covered by Peng's BSDE
method and requires a quite new approach. The uniqueness of the viscosity
solution will be shown in the next section for the class of bounded
continuous functions. We first recall the definition of a viscosity solution
of (\ref{HJBI}). The interested reader is referred to Crandall, Ishii, and
Lions \cite{CIL}.

\begin{definition}
\label{dvis1} Let $u\left( t,x\right) \in C\left( \left[ 0,T\right] \times
\mathbb{R}^{n}\right) $ and $\left( t,x\right) \in \left[ 0,T\right] \times
\mathbb{R}^{n}.$ For every $\varphi \in C^{1,2}\left( \left[ 0,T\right]
\times \mathbb{R}^{n}\right) $

\noindent (1) for each local maximum point $\left( t_{0},x_{0}\right) $ of $%
u-\varphi $ in the interior of $\left[ 0,T\right] \times \mathbb{R}^{n},$ we
have%
\begin{equation}
\max \left\{ V-\mathcal{H}_{\inf }^{\chi }V,\min \left[ -\frac{\partial }{%
\partial t}\varphi -H\left( t,x,\varphi ,D\varphi ,D^{2}\varphi \right) ,V-%
\mathcal{H}_{\sup }^{c}V\right] \right\} \leq 0  \label{vis1}
\end{equation}%
and for each $x\in \mathbb{R}^{n},$ we have%
\begin{equation}
\max \left\{ V\left( T,x\right) -\mathcal{H}_{\inf }^{\chi }V\left(
T,x\right) ,\min \left[ V\left( T,x\right) -\Phi \left( x\right) ,V\left(
T,x\right) -\mathcal{H}_{\inf }^{c}V\left( T,x\right) \right] \right\} \leq 0
\label{TT1}
\end{equation}%
i.e., $u$ is a subsolution to HJBI equation \emph{(\ref{HJBI})}.

\noindent (2) for each local minimum point $\left( t_{0},x_{0}\right) $ of $%
u-\varphi $ in the interior of $\left[ 0,T\right] \times \mathbb{R}^{n},$ we
have%
\begin{equation}
\max \left\{ V-\mathcal{H}_{\inf }^{\chi }V,\min \left[ -\frac{\partial }{%
\partial t}\varphi -H\left( t,x,\varphi ,D\varphi ,D^{2}\varphi \right) ,V-%
\mathcal{H}_{\sup }^{c}V\right] \right\} \geq 0  \label{vis2}
\end{equation}%
and for each $x\in \mathbb{R}^{n},$ we have%
\begin{equation}
\max \left\{ V\left( T,x\right) -\mathcal{H}_{\inf }^{\chi }V\left(
T,x\right) ,\min \left[ V\left( T,x\right) -\Phi \left( x\right) ,V\left(
T,x\right) -\mathcal{H}_{\inf }^{c}V\left( T,x\right) \right] \right\} \geq 0
\label{TT2}
\end{equation}%
i.e., $u$ is a supersolution to HJBI equation \emph{(\ref{HJBI})}.

\noindent (3) $u\left( t,x\right) \in C\left( \left[ 0,T\right] \times
\mathbb{R}^{n}\right) $ is said to be a viscosity solution of \emph{(\ref%
{HJBI})} if it is both a viscosity sub and supersolution.
\end{definition}

We have the other definition which will be useful to verify the viscosity
solutions.

\begin{definition}
Let $u\left( t,x\right) \in C\left( \left[ 0,T\right] \times \mathbb{R}%
^{n}\right) $ and $\left( t,x\right) \in \left[ 0,T\right] \times \mathbb{R}%
^{n}$. We denote by $\mathcal{P}^{2,+}u\left( t,x\right) $, the
\textquotedblleft parabolic superjet\textquotedblright\ of $u$ at $\left(
t,x\right) $ the set of triples $\left( p,q,X\right) \in \mathbb{R}\times
\mathbb{R}^{n}\times \mathbb{S}^{n}$ which are such that%
\begin{eqnarray*}
u\left( s,y\right) &\leq &u\left( t,x\right) +p\left( s-t\right)
+\left\langle q,x-y\right\rangle \\
&&+\frac{1}{2}\left\langle X\left( y-x\right) ,y-x\right\rangle +o\left(
\left\vert s-t\right\vert +\left\vert y-x\right\vert ^{2}\right) .
\end{eqnarray*}%
Similarly, we denote by $\mathcal{P}^{2,-}u\left( t,x\right) ,$ the
\textquotedblright parabolic subjet\textquotedblright\ of $u$ at $\left(
t,x\right) $ the set of triples $\left( p,q,X\right) \in \mathbb{R}\times
\mathbb{R}^{n}\times \mathbb{S}^{n}$ which are such that%
\begin{eqnarray*}
u\left( s,y\right) &\geq &u\left( t,x\right) +p\left( s-t\right)
+\left\langle q,x-y\right\rangle \\
&&+\frac{1}{2}\left\langle X\left( y-x\right) ,y-x\right\rangle +o\left(
\left\vert s-t\right\vert +\left\vert y-x\right\vert ^{2}\right) .
\end{eqnarray*}
\end{definition}

\begin{definition}
\label{dvis2} (i) It can be said $V\left( t,x\right) \in C\left( \left[ 0,T%
\right] \times \mathbb{R}^{n}\right) $ is a viscosity subsolution of \emph{(%
\ref{HJBI}) }if\emph{\ }at any point $\left( t,x\right) \in \left[ 0,T\right]
\times \mathbb{R}^{n}$, for any $\left( p,q,X\right) \in \mathcal{P}%
^{2,+}V\left( t,x\right) $,%
\begin{equation}
\max \left\{ V-\mathcal{H}_{\inf }^{\chi }V,\min \left[ -p-H\left(
t,x,V\left( t,x\right) ,q,X\right) ,V-\mathcal{H}_{\sup }^{c}V\right]
\right\} \leq 0  \label{other1}
\end{equation}%
and for each $x\in \mathbb{R}^{n},$ it holds%
\begin{equation}
\max \left\{ V\left( T,x\right) -\mathcal{H}_{\inf }^{\chi }V\left(
T,x\right) ,\min \left[ V\left( T,x\right) -\Phi \left( x\right) ,V\left(
T,x\right) -\mathcal{H}_{\inf }^{c}V\left( T,x\right) \right] \right\} \leq
0.  \label{T1}
\end{equation}

\noindent (ii) It can be said $V\left( t,x\right) \in C\left( \left[ 0,T%
\right] \times \mathbb{R}^{n}\right) $ is a viscosity supersolution of \emph{%
(\ref{HJBI})} if at any point $\left( t,x\right) \in \left[ 0,T\right]
\times \mathbb{R}^{n}$, for any $\left( p,q,X\right) \in \mathcal{P}%
^{2,+}V\left( t,x\right) $,%
\begin{equation}
\max \left\{ V-\mathcal{H}_{\inf }^{\chi }V,\min \left[ -p-H\left(
t,x,V\left( t,x\right) ,q,X\right) ,V-\mathcal{H}_{\sup }^{c}V\right]
\right\} \geq 0  \label{other2}
\end{equation}%
and for each $x\in \mathbb{R}^{n},$ we have%
\begin{equation}
\max \left\{ V\left( T,x\right) -\mathcal{H}_{\inf }^{\chi }V\left(
T,x\right) ,\min \left[ V\left( T,x\right) -\Phi \left( x\right) ,V\left(
T,x\right) -\mathcal{H}_{\inf }^{c}V\left( T,x\right) \right] \right\} \geq
0.  \label{T2}
\end{equation}%
\noindent (iii) It can be said $u\left( t,x\right) \in C\left( \left[ 0,T%
\right] \times \mathbb{R}^{n}\right) $ is a viscosity solution of \emph{(\ref%
{HJBI})} if it is both a viscosity sub and super solution.
\end{definition}

\begin{remark}
Definition \ref{dvis1} and \ref{dvis2} are equivalent to each other. For
more details, see Fleming and Soner \cite{FS}, Lemma 4.1 (page 211).
\end{remark}

We now introduce the lower and upper obstacles with the help of the
following lemms.

\begin{lemma}
\label{obs}Assume $\emph{(A1)}$\emph{-}$\emph{(A3)}$ are in force. Given any
$\left( t,x\right) \in \left( 0,T\right] \times \mathbb{R}^{n},$ the lower
and upper value functions satisfy the following equation:%
\begin{equation*}
\max \left\{ \min \left[ 0,V\left( t,x\right) -\mathcal{H}_{\sup
}^{c}V\left( t,x\right) \right] ,V\left( t,x\right) -\mathcal{H}_{\inf
}^{\chi }V\left( t,x\right) \right\} =0.
\end{equation*}
\end{lemma}

\paragraph{Proof}

From Lemma \ref{deter}, (\ref{st0}) can be expressed as
\begin{eqnarray*}
V^{-}\left( t,x\right) &=&\inf_{\rho \in \mathcal{T}_{t,+\infty }\eta \in
\mathcal{F}_{\rho }}\sup_{\tau \in \mathcal{T}_{t,+\infty },\xi \in \mathcal{%
F}_{\tau }}\mathbb{E}\Big [-c\left( t,\xi \right) \mathbf{1}_{\left\{ \tau
=t\right\} }\mathbf{1}_{\left\{ \rho =+\infty \right\} } \\
&&+\chi \left( t,\eta \right) \mathbf{1}_{\left\{ \rho =t\right\}
}+V^{-}\left( t,X_{t}^{t,x;\xi \mathbf{1}_{\left[ \tau ,T\right] },\eta
\mathbf{1}_{\left[ \rho ,T\right] }}\right) \Big ].
\end{eqnarray*}%
The remainder of the proof is the same as Lemma 5.3 from \cite{co}. We omit
it.\hfill $\Box $

\begin{remark}
We have $V^{-}\left( t,x\right) \leq \mathcal{H}_{\inf }^{\chi }V\left(
t,x\right) $ on $\left( 0,T\right] \times \mathbb{R}^{n}$ from Lemma \ref%
{obs}. Besides, whenever $V\left( t,x\right) \leq \mathcal{H}_{\inf }^{\chi
}V\left( t,x\right) $, then $\mathcal{H}_{\sup }^{c}V\left( t,x\right) \leq
V^{-}\left( t,x\right) $ and $\mathcal{H}_{\sup }^{c}V\left( t,x\right) \leq
\mathcal{H}_{\inf }^{\chi }V\left( t,x\right) $. So we may regard $\mathcal{H%
}_{\sup }^{c}V\left( t,x\right) $ as a lower obstacle and $\mathcal{H}_{\inf
}^{\chi }V\left( t,x\right) $ as an upper obstacle. Both of them are
implicit forms, since they depend on $V^{-}$. The same remark applies to $%
V^{+}$ likewise.
\end{remark}

We shall prove that the two value functions satisfy, in the viscosity sense,
the terminal condition.

\begin{lemma}
\label{tervis}Assume assumptions $\emph{(A1)}$\emph{-}$\emph{(A3)}$ are in
force. The lower value function $V^{-}\left( T,x\right) $ is a viscosity
solution of \emph{(\ref{HJBI}).}
\end{lemma}

\paragraph{Proof}

We shall prove%
\begin{equation*}
\max \left\{ V\left( T,x\right) -\mathcal{H}_{\inf }^{\chi }V\left(
T,x\right) ,\min \left[ V\left( T,x\right) -\Phi \left( x\right) ,V\left(
T,x\right) -\mathcal{H}_{\inf }^{c}V\left( T,x\right) \right] \right\} \geq
0.
\end{equation*}%
From Lemma \ref{ter}, we have%
\begin{eqnarray*}
V^{-}\left( t,x\right) &=&\inf_{\rho \in \mathcal{T}_{t,+\infty }\eta \in
\mathcal{F}_{\rho }}\sup_{\tau \in \mathcal{T}_{t,+\infty },\xi \in \mathcal{%
F}_{\tau }}\mathbb{E}\Big [\Big (-c\left( t,\xi \right) \mathbf{1}_{\left\{
\tau =t\right\} }\mathbf{1}_{\left\{ \rho =+\infty \right\} } \\
&&+\chi \left( t,\eta \right) \mathbf{1}_{\left\{ \rho =t\right\}
}+V^{-}\left( t,X_{t}^{t,x;\xi \mathbf{1}_{\left[ \tau ,T\right] },\eta
\mathbf{1}_{\left[ \rho ,T\right] }}\right) \Big )\left( 1-\mathbf{1}%
_{\left\{ \tau =+\infty ,\rho =+\infty \right\} }\right) \\
&&+\left( \int_{t}^{T}f\left(
s,X_{s}^{t,x;u_{0},v_{0}},Y_{s}^{t,x;u_{0},v_{0}},Z_{s}^{t,x;u_{0},v_{0}}%
\right) \mathrm{d}s+\Phi \left( X_{T}^{t,x;u_{0},v_{0}}\right) \right)
\mathbf{1}_{\left\{ \tau =+\infty ,\rho =+\infty \right\} }\Big ].
\end{eqnarray*}%
Thanks to (A1)-(A2), it follows that%
\begin{eqnarray*}
&&\mathbb{E}\left[ \left\vert \int_{t}^{T}f\left(
s,X_{s}^{t,x;u_{0},v_{0}},Y_{s}^{t,x;u_{0},v_{0}},Z_{s}^{t,x;u_{0},v_{0}}%
\right) \mathrm{d}s\right\vert \right] \\
&\leq &\left( T-t\right) ^{\frac{1}{2}}\mathbb{E}\left[ \int_{t}^{T}\left%
\vert f\left(
s,X_{s}^{t,x;u_{0},v_{0}},Y_{s}^{t,x;u_{0},v_{0}},Z_{s}^{t,x;u_{0},v_{0}}%
\right) \right\vert ^{2}\mathrm{d}s\right] ^{\frac{1}{2}} \\
&\leq &C\left( T-t\right) ^{\frac{1}{2}},
\end{eqnarray*}%
and
\begin{eqnarray*}
\mathbb{E}\left[ \left\vert \Phi \left( X_{T}^{t,x;u_{0},v_{0}}\right) -\Phi
\left( x\right) \right\vert \right] &\leq &C\mathbb{E}\left[ \left\vert
X_{T}^{t,x;u_{0},v_{0}}-x\right\vert ^{2}\right] ^{\frac{1}{2}} \\
&\leq &C\left( T-t\right) ^{\frac{1}{2}},\text{ uniformly in }u_{0},v_{0}.
\end{eqnarray*}%
Therefore,
\begin{eqnarray*}
V^{-}\left( t,x\right) &\geq &-\inf_{\rho \in \mathcal{T}_{t,+\infty }\eta
\in \mathcal{F}_{\rho }}\sup_{\tau \in \mathcal{T}_{t,+\infty },\xi \in
\mathcal{F}_{\tau }}\mathbb{E}\Big [\Big (-c\left( t,\xi \right) \mathbf{1}%
_{\left\{ \tau =t\right\} }\mathbf{1}_{\left\{ \rho =+\infty \right\} } \\
&&+\chi \left( t,\eta \right) \mathbf{1}_{\left\{ \rho =t\right\}
}+V^{-}\left( t,X_{t}^{t,x;\xi \mathbf{1}_{\left[ \tau ,T\right] },\eta
\mathbf{1}_{\left[ \rho ,T\right] }}\right) \Big )\left( 1-\mathbf{1}%
_{\left\{ \tau =+\infty ,\rho =+\infty \right\} }\right) \\
&&+\Phi \left( x\right) \mathbf{1}_{\left\{ \tau =+\infty ,\rho =+\infty
\right\} }\Big ]-C\left( T-t\right) ^{\frac{1}{2}}.
\end{eqnarray*}%
Repeating the method in Lemma \ref{obs}, we have
\begin{eqnarray*}
&&\max \left\{ \min \left[ V^{-}\left( t,x\right) -\Phi \left( x\right)
,V^{-}\left( t,x\right) -\mathcal{H}_{\sup }^{c}V^{-}\left( t,x\right) %
\right] ,V^{-}\left( t,x\right) -\mathcal{H}_{\inf }^{\chi }V^{-}\left(
t,x\right) \right\} \\
&\geq &-C\left( T-t\right) ^{\frac{1}{2}}.
\end{eqnarray*}%
According to (A3), namely the $1/2$-H\"{o}lder continuity in time for $c,$ $%
\chi $ and $V^{-},$ we deduce that
\begin{eqnarray}
&&\max \left\{ V\left( T,x\right) -\mathcal{H}_{\inf }^{\chi }V\left(
T,x\right) ,\min \left[ V\left( T,x\right) -\Phi \left( x\right) ,V\left(
T,x\right) -\mathcal{H}_{\inf }^{c}V\left( T,x\right) \right] \right\}
\notag \\
&&+C\left( T-t\right) ^{\frac{1}{2}}  \notag \\
&\geq &\max \left\{ \min \left[ V^{-}\left( t,x\right) -\Phi \left( x\right)
,V^{-}\left( t,x\right) -\mathcal{H}_{\sup }^{c}V^{-}\left( t,x\right) %
\right] ,V^{-}\left( t,x\right) -\mathcal{H}_{\inf }^{\chi }V^{-}\left(
t,x\right) \right\}  \notag \\
&\geq &-C_{1}\left( T-t\right) ^{\frac{1}{2}}.  \label{tervis1}
\end{eqnarray}%
for some $C_{1}>0.$ Then, letting $t=T$ in (\ref{tervis1}) ends the
proof.\hfill $\Box $

We first prove that the lower value function $V^{-}\left( t,x\right) $ is a
viscosity solution of (\ref{HJBI}).

\begin{theorem}
Assume assumptions $\emph{(A1)}$\emph{-}$\emph{(A3)}$ are in force, the
lower value function $V^{-}\left( t,x\right) $ is a viscosity solution of
\emph{(\ref{HJBI})}.
\end{theorem}

\paragraph{Proof}

We first show that the lower value function $V^{-}$ is a viscosity solution
to (\ref{HJBI}); the other case is analogous.

In Lemma \ref{tervis}, we have proved that $V^{-}$ satisfies, in the
viscosity sense, the terminal condition, namely (\ref{T1}) and (\ref{T1}).
Therefore, we have only to address (\ref{other1}). From Proposition \ref%
{vlip}, $V^{-}$ is continuous on $\left[ 0,T\right) \times \mathbb{R}^{n}$.
Thus we begin by proving that $V^{-}$ is a viscosity supersolution. By
virtue of Lemma \ref{obs}, we have to show that, given $\left( \bar{t},\bar{x%
}\right) \in \left[ 0,T\right) \times \mathbb{R}^{n}$ such that $\mathcal{H}%
_{\sup }^{c}V\left( \bar{t},\bar{x}\right) \leq V^{-}\left( \bar{t},\bar{x}%
\right) $ and $V^{-}\left( \bar{t},\bar{x}\right) \leq \mathcal{H}_{\inf
}^{\chi }V\left( \bar{t},\bar{x}\right) ,$ then for every $\varphi \in
C^{1,2}\left( \left[ 0,T\right) \times \mathbb{R}^{n}\right) $, such that $%
\left( \bar{t},\bar{x}\right) $ is a local minimum of $V^{-}-\varphi $, we
have%
\begin{equation}
\frac{\partial }{\partial t}\varphi \left( \bar{t},\bar{x}\right) +H\left(
t,x,\varphi \left( \bar{t},\bar{x}\right) ,D\varphi \left( \bar{t},\bar{x}%
\right) ,D^{2}\varphi \left( \bar{t},\bar{x}\right) \right) \leq 0,
\label{supvis}
\end{equation}%
where $H$ is defiend in (\ref{H}).

Without loss of generality, postulate $V^{-}\left( \bar{t},\bar{x}\right)
=\varphi \left( \bar{t},\bar{x}\right) .$ Let
\begin{equation}
\lambda +V^{-}\left( \bar{t},\bar{x}\right) =\mathcal{H}_{\inf }^{\chi
}V\left( \bar{t},\bar{x}\right) =\inf_{y\in \mathcal{V}}\left[ V^{-}\left(
\bar{t},\bar{x}+y\right) +\chi \left( \bar{t},y\right) \right] .
\label{supesti}
\end{equation}%
We proceed as in \cite{co} to derive the following result: For every random
variable $\eta $, $\mathcal{F}_{s}$-measurable and values in $\mathcal{V}$,
there exists $C>0,$%
\begin{equation*}
\mathbb{E}^{\mathcal{F}_{s}}\left[ V\left( s,X_{s}^{\bar{t},\bar{x}}\right) %
\right] \leq \mathbb{E}^{\mathcal{F}_{s}}\left[ V\left( s,X_{s}^{\bar{t},%
\bar{x}}+\eta \right) +\chi \left( s,\eta \right) \right] +C\left\vert s-%
\bar{t}\right\vert ^{\frac{1}{2}}-\lambda
\end{equation*}%
with $X_{s}^{\bar{t},\bar{x}}=X_{s}^{\bar{t},\bar{x},u_{0},v_{0}}$ for all $%
s\in \left[ \bar{t},T\right] $, $P$-a.s., where $u_{0}$ and $v_{0}$ denote
the controls with no impulses.

Next recall%
\begin{equation*}
V^{-}\left( \bar{t},\bar{x}\right) =\inf_{\beta \in \mathcal{B}_{\bar{t}%
,T}}\sup_{u\in \mathcal{U}_{\bar{t},T}}G_{\bar{t},\bar{t}+\delta }^{\bar{t},%
\bar{x};u,\beta \left( u\right) }\left[ V^{-}\left( \bar{t}+\delta ,X_{\bar{t%
}+\delta }^{\bar{t},\bar{x};u,\beta \left( u\right) }\right) +\Theta _{\bar{t%
}+\delta }^{u,\beta \left( u\right) }\right] ,
\end{equation*}%
where
\begin{equation*}
\Theta _{\bar{t}+\delta }^{u,v}:=\sum_{l\geq 1}\chi \left( \rho _{l},\eta
_{l}\right) \mathbf{1}_{\left\{ \rho _{l}\leq \bar{t}+\delta \right\}
}-\sum_{m\geq 1}c\left( \tau _{m},\xi _{m}\right) \mathbf{1}_{\left\{ \tau
_{m}\leq \bar{t}+\delta \right\} }\prod_{l\geq 1}\mathbf{1}_{\left\{ \tau
_{m}\neq \rho _{l}\right\} }.
\end{equation*}%
From the definition of $V^{-}\left( \bar{t},\bar{x}\right) ,$ we have
\begin{eqnarray*}
V^{-}\left( \bar{t},\bar{x}\right) &=&\inf_{\beta \in \mathcal{B}_{\bar{t}%
,T}}\sup_{u\in \mathcal{U}_{\bar{t},T}}G_{\bar{t},\bar{t}+\delta }^{\bar{t},%
\bar{x};u,\beta \left( u\right) }\left[ V^{-}\left( \bar{t}+\delta ,X_{\bar{t%
}+\delta }^{\bar{t},\bar{x};u,\beta \left( u\right) }\right) +\Theta _{\bar{t%
}+\delta }^{u,\beta \left( u\right) }\right] \\
&\geq &G_{\bar{t},\bar{t}+\delta }^{\bar{t},\bar{x};u_{0},\beta
^{\varepsilon }\left( u_{0}\right) }\left[ V^{-}\left( \bar{t}+\delta ,X_{%
\bar{t}+\delta }^{\bar{t},\bar{x};u_{0},\beta ^{\varepsilon }\left(
u_{0}\right) }\right) +\sum_{l\geq 1}\chi \left( \rho _{l},\eta _{l}\right)
\mathbf{1}_{\left\{ \rho _{l}\leq \bar{t}+\delta \right\} }\right]
-\varepsilon
\end{eqnarray*}%
for some $\beta ^{\varepsilon }\in \mathcal{B}_{\bar{t},T},$ with $%
v^{\varepsilon }:=\beta ^{\varepsilon }\left( u_{0}\right) =\sum_{l\geq
1}\eta _{l}^{\varepsilon }\mathbf{1}_{\left[ \rho _{l}^{\varepsilon },T%
\right] }\in \mathcal{V}_{\bar{t},T}.$ Note that
\begin{equation*}
\sum_{l\geq 1}\chi \left( \rho _{l},\eta _{l}\right) \mathbf{1}_{\left\{
\rho _{l}\leq \bar{t}+\delta \right\} }=\sum_{l\geq 1}^{\mu _{\left[ \bar{t},%
\bar{t}+\delta \right] }}\chi \left( \rho _{l},\eta _{l}\right) \mathbf{.}
\end{equation*}%
By (A1)-(A2), we have the following estimate%
\begin{eqnarray*}
&&\mathbb{E}^{\mathcal{F}_{\bar{t}+\delta }}\left[ \left\vert V^{-}\left(
\bar{t}+\delta ,X_{\bar{t}+\delta }^{\bar{t},\bar{x};u_{0},v^{\varepsilon
}}\right) -V^{-}\left( \bar{t}+\delta ,X_{s}^{\bar{t},\bar{x}}+\sum_{l\geq
1}^{\mu _{\left[ \bar{t},\bar{t}+\delta \right] }}\chi \left( \rho _{l},\eta
_{l}\right) \right) \right\vert \right] \\
&\leq &C\delta ^{\frac{1}{2}}\mathbb{E}^{\mathcal{F}_{\bar{t}+\delta }}\left[
\mathbf{1}_{\left\{ \mu _{\left[ \bar{t},\bar{t}+\delta \right] }\geq
1\right\} }\right] .
\end{eqnarray*}%
As a consequence, using (A3) and (\ref{supvis}), we deduce
\begin{eqnarray*}
&&\mathbb{E}^{\mathcal{F}_{\bar{t}+\delta }}\left[ V^{-}\left( \bar{t}%
+\delta ,X_{\bar{t}+\delta }^{\bar{t},\bar{x}}+\sum_{l\geq 1}^{\mu _{\left[
\bar{t},\bar{t}+\delta \right] }}\chi \left( \rho _{l},\eta _{l}\right)
\right) +\sum_{l\geq 1}^{\mu _{\left[ \bar{t},\bar{t}+\delta \right] }}\chi
\left( \rho _{l},\eta _{l}\right) \right] \\
&\geq &\mathbb{E}^{\mathcal{F}_{\bar{t}+\delta }}\left[ V^{-}\left( \bar{t}%
+\delta ,X_{\bar{t}+\delta }^{\bar{t},\bar{x}}\right) +\left( \lambda
-C\delta ^{\frac{1}{2}}\right) \mathbf{1}_{\left\{ \mu _{\left[ \bar{t},\bar{%
t}+\delta \right] }\geq 1\right\} }\right] .
\end{eqnarray*}%
Therefore, applying comparison theorem (Proposition 2.6 in \cite{BBP}), we
find%
\begin{eqnarray*}
&&G_{\bar{t},\bar{t}+\delta }^{\bar{t},\bar{x};u_{0},v^{\varepsilon }}\left[
V^{-}\left( \bar{t}+\delta ,X_{\bar{t}+\delta }^{\bar{t},\bar{x}%
;u_{0},v^{\varepsilon }}\right) +\sum_{l\geq 1}\chi \left( \rho _{l},\eta
_{l}\right) \mathbf{1}_{\left\{ \rho _{l}\leq \bar{t}+\delta \right\} }%
\right] \\
&\geq &G_{\bar{t},\bar{t}+\delta }^{\bar{t},\bar{x};u_{0},v^{\varepsilon }}%
\left[ V^{-}\left( \bar{t}+\delta ,X_{\bar{t}+\delta }^{\bar{t},\bar{x}%
}\right) +\left( \lambda -C\delta ^{\frac{1}{2}}\right) \mathbf{1}_{\left\{
\mu _{\left[ \bar{t},\bar{t}+\delta \right] }\geq 1\right\} }\right] .
\end{eqnarray*}%
Thus
\begin{equation*}
V^{-}\left( \bar{t},\bar{x}\right) \geq G_{\bar{t},\bar{t}+\delta }^{\bar{t},%
\bar{x};u_{0},v^{\varepsilon }}\left[ V^{-}\left( \bar{t}+\delta ,X_{\bar{t}%
+\delta }^{\bar{t},\bar{x}}\right) +\left( \lambda -C\delta ^{\frac{1}{2}%
}\right) \mathbf{1}_{\left\{ \mu _{\left[ \bar{t},\bar{t}+\delta \right]
}\geq 1\right\} }\right] -\varepsilon .
\end{equation*}%
From the boundedness of $f$ we deduce%
\begin{eqnarray*}
&&f\left( s,X_{s}^{\bar{t},\bar{x};u_{0},v^{\varepsilon }},y,z\right) \\
&=&f\left( s,X_{s}^{\bar{t},\bar{x};u_{0},v^{\varepsilon }},y,z\right)
-f\left( s,X_{s}^{\bar{t},\bar{x}},y,z\right) +f\left( s,X_{s}^{\bar{t},\bar{%
x}},y,z\right) \\
&\geq &f\left( s,X_{s}^{\bar{t},\bar{x}},y,z\right) -C,\text{ for }\left(
y,z\right) \in \mathbb{R}\times \mathbb{R}^{d}.
\end{eqnarray*}%
Applying comparison theorem (Proposition 2.6 in \cite{BBP}) again, we have $%
V^{-}\left( \bar{t},\bar{x}\right) \geq \mathcal{Y}_{\bar{t}},$ where $%
\mathcal{Y}_{\bar{t}}$ is the solution to the following BSDE:%
\begin{eqnarray*}
\mathcal{Y}_{\bar{t}} &=&V^{-}\left( \bar{t}+\delta ,X_{\bar{t}+\delta }^{%
\bar{t},\bar{x}}\right) +\left( \lambda -C\delta ^{\frac{1}{2}}-C\delta
\right) \mathbf{1}_{\left\{ \mu _{\left[ \bar{t},\bar{t}+\delta \right]
}\geq 1\right\} } \\
&&+\int_{\bar{t}}^{\bar{t}+\delta }f\left( s,X_{s}^{\bar{t},\bar{x}},%
\mathcal{Y}_{s},\mathcal{Z}_{s}\right) \mathrm{d}s-\int_{\bar{t}}^{\bar{t}%
+\delta }\mathcal{Z}_{s}\mathrm{d}W_{s}.
\end{eqnarray*}%
We shall take $\delta $ sufficiently small. Indeed, there exists $\bar{\delta%
}>0$ such that for $\zeta \in \left( 0,\bar{\delta}\right) ,$ we have $%
\lambda -C\delta ^{\frac{1}{2}}-C\delta \geq 0.$ Immediately, by Proposition
2.6 in \cite{BBP}, it follows
\begin{equation}
V^{-}\left( \bar{t},\bar{x}\right) \geq G_{\bar{t},\bar{t}+\zeta }^{\bar{t},%
\bar{x}}\left[ V^{-}\left( \bar{t}+\zeta ,X_{\bar{t}+\zeta }^{\bar{t},\bar{x}%
}\right) \right] -\varepsilon .  \label{key}
\end{equation}%
To abbreviate notations we set, for some arbitrarily chosen but fixed $%
\varphi \in C^{1,2}\left( \left[ 0,T\right) \times \mathbb{R}^{n}\right) ,$%
\begin{eqnarray*}
F\left( s,x,y,z\right) &=&\frac{\partial }{\partial s}\varphi \left(
s,x\right) +\frac{1}{2}\text{Tr}\left( \sigma \sigma ^{\top }\left(
s,x\right) D^{2}\varphi \right) +\left\langle D\varphi ,b\left( s,x\right)
\right\rangle \\
&&+f\left( s,x,y+\varphi \left( s,x\right) ,z+D\varphi \left( s,x\right)
\cdot \sigma \left( s,x\right) \right) ,
\end{eqnarray*}%
for $\left( s,x,y,z\right) \in \left[ 0,T\right] \times \mathbb{R}^{n}%
\mathbb{\times R\times R}^{d}.$

Let us consider the following BSDE:%
\begin{equation}
\left\{
\begin{array}{rcl}
-\mathrm{d}Y_{s}^{1} & = & F\left( s,X_{s}^{\bar{t},\bar{x}%
},Y_{s}^{1},Z_{s}^{1}\right) \mathrm{d}s-Z_{s}^{1}\mathrm{d}W_{s}, \\
Y_{\bar{t}+\zeta }^{1} & = & 0.%
\end{array}%
\right.  \label{vbsde1}
\end{equation}%
It is not hard to check that $F\left( s,X_{s}^{\bar{t},\bar{x}},y,z\right) $
satisfies (A1) and (A2). Thus, BSDE (\ref{vbsde1}) admits a unique adapted
strong solution. We can characterize the solution process $Y_{s}^{1}$ as
follows.%
\begin{equation}
Y_{s}^{1}=G_{s,\bar{t}+\zeta }^{\bar{t},\bar{x}}\left[ \varphi \left( \bar{t}%
+\zeta ,X_{\bar{t}+\zeta }^{\bar{t},\bar{x}}\right) \right] -\varphi \left(
s,X_{s}^{\bar{t},\bar{x}}\right) .  \label{key2}
\end{equation}%
Indeed, $G_{s,\bar{t}+\zeta }^{\bar{t},\bar{x}}\left[ \varphi \left( \bar{t}%
+\zeta ,X_{\bar{t}+\zeta }^{\bar{t},\bar{x}}\right) \right] $ is defined by
the solution of the following BSDE:%
\begin{equation}
\left\{
\begin{array}{rcl}
-\mathrm{d}\bar{Y}_{s} & = & f\left( s,X_{s}^{\bar{t},\bar{x}},\bar{Y}_{s},%
\bar{Z}_{s}\right) \mathrm{d}s-\bar{Z}_{s}\mathrm{d}W_{s}, \\
\bar{Y}_{\bar{t}+\zeta } & = & \varphi \left( \bar{t}+\zeta ,X_{\bar{t}%
+\zeta }^{\bar{t},\bar{x}}\right) .%
\end{array}%
\right.  \label{ebsde}
\end{equation}%
Therefore, one just need to prove $\bar{Y}_{s}-\varphi \left( s,X_{s}^{\bar{t%
},\bar{x}}\right) =Y_{s}^{1}$. Applying It\^{o}'s formula to $\varphi \left(
s,X_{s}^{\bar{t},\bar{x}}\right) ,$ we obtain $\mathrm{d}\left[ \bar{Y}%
_{s}-\varphi \left( s,X_{s}^{\bar{t},\bar{x}}\right) \right] =\mathrm{d}%
Y_{s}^{1}$, and at the terminal time $\bar{Y}_{\bar{t}+\zeta }-\varphi
\left( \bar{t}+\zeta ,X_{\bar{t}+\zeta }^{\bar{t},\bar{x}}\right) =Y_{\bar{t}%
+\zeta }^{1}=0$, as a result, they are equal in the interval $\left[ \bar{t},%
\bar{t}+\zeta \right] $.

Now let us introduce a more simpler BSDE than (\ref{vbsde1}), i.e., $X_{s}^{%
\bar{t},\bar{x}}$ of the equation (\ref{vbsde1}) is taken place by $x$:%
\begin{equation}
\left\{
\begin{array}{rcl}
-\mathrm{d}Y_{s}^{2} & = & F\left( s,\bar{x},Y_{s}^{2},Z_{s}^{2}\right)
\mathrm{d}s-Z_{s}^{2}\mathrm{d}W_{s}, \\
Y_{\bar{t}+\zeta }^{2} & = & 0.%
\end{array}%
\right.  \label{vbsde2}
\end{equation}%
Notice that $F$ is a deterministic function of $(s,x,y,z)$ therefore $\left(
Y_{s}^{2},Z_{s}^{2}\right) =\left( Y_{0}\left( s\right) ,0\right) $ where $%
Y_{0}\left( s\right) $ is the solution of the ODE:%
\begin{equation}
\left\{
\begin{array}{rcl}
-\dot{Y}_{0}\left( s\right) & = & F\left( s,\bar{x},Y_{0}\left( s\right)
,0\right) \mathrm{d}s,\text{ }s\in \left[ \bar{t},\bar{t}+\zeta \right] , \\
Y_{0}\left( \bar{t}+\zeta \right) & = & 0.%
\end{array}%
\right.  \label{ode}
\end{equation}%
The following result indicates that the difference of the solutions of (\ref%
{vbsde1}) and (\ref{vbsde2}) can be neglected whenever $\zeta $ is
sufficiently small enough. From the classical estimate on SDE, we have $%
\mathbb{E}\left[ \sup_{s\in \left[ \bar{t},\bar{t}+\zeta \right] }\left\vert
X_{s}^{\bar{t},\bar{x}}\right\vert ^{p}\right] \leq C\left( 1+\left\vert
\bar{x}\right\vert ^{p}\right) .$ Moreover, applying B-D-G inequality, we
get $\mathbb{E}\left[ \sup_{s\in \left[ \bar{t},\bar{t}+\zeta \right]
}\left\vert X_{s}^{\bar{t},\bar{x}}-\bar{x}\right\vert ^{2}\right] \leq
C\zeta .$ Hence, when $\zeta \rightarrow 0,$ the following random variable $%
\kappa ^{\zeta }:=\sup_{s\in \left[ \bar{t},\bar{t}+\zeta \right]
}\left\vert X_{s}^{\bar{t},\bar{x}}-\bar{x}\right\vert $ converges monotone
to $0$. On the one hand, employing Proposition 3.2 in \cite{BDHPS} to BSDEs (%
\ref{vbsde1}) and(\ref{vbsde2}), we have%
\begin{equation*}
\mathbb{E}\left[ \int_{\bar{t}}^{\bar{t}+\zeta }\left\vert
Y_{s}^{1}-Y_{s}^{2}\right\vert ^{2}+\left\vert Z_{s}^{1}\right\vert ^{2}%
\right] \mathrm{d}s\leq C\mathbb{E}\left[ \int_{\bar{t}}^{\bar{t}+\zeta
}\varpi \left( \left\vert X_{s}^{\bar{t},\bar{x}}-x\right\vert \right) ^{2}%
\mathrm{d}s\right] \leq C\zeta \mathbb{E}\varpi \left( \kappa ^{\zeta
}\right) ^{2}.
\end{equation*}%
On the other hand, from Lemma \ref{deter}, we have
\begin{eqnarray}
\left\vert Y_{\bar{t}}^{1}-Y_{\bar{t}}^{2}\right\vert &=&\left\vert \mathbb{E%
}\left\vert Y_{\bar{t}}^{1}\right\vert \right\vert  \notag \\
&=&\left\vert \mathbb{E}\left[ \int_{\bar{t}}^{\bar{t}+\zeta }\left[ F\left(
s,X_{s}^{\bar{t},\bar{x}},Y_{s}^{1},Z_{s}^{1}\right) -F\left( s,\bar{x}%
,Y_{s}^{2}\right) \right] \mathrm{d}s\right] \right\vert  \notag \\
&\leq &C\mathbb{E}\left[ \int_{\bar{t}}^{\bar{t}+\zeta }\left( \varpi \left(
\left\vert X_{s}^{\bar{t},\bar{x}}-x\right\vert \right) +\left\vert
Y_{s}^{1}-Y_{s}^{2}\right\vert +\left\vert Z_{s}^{1}\right\vert \right)
\mathrm{d}s\right]  \notag \\
&\leq &C\zeta \mathbb{E}\varpi \left( \kappa ^{\zeta }\right) ^{2}+C\zeta ^{%
\frac{1}{2}}\left\{ \mathbb{E}\left[ \int_{\bar{t}}^{\bar{t}+\zeta }\left(
\left\vert Y_{s}^{1}-Y_{s}^{2}\right\vert ^{2}+\left\vert
Z_{s}^{1}\right\vert ^{2}\right) \mathrm{d}s\right] \right\} ^{\frac{1}{2}}
\notag \\
&\leq &C\zeta \mathbb{E}\left[ \varpi \left( \kappa ^{\zeta }\right)
^{2}+\varpi \left( \kappa ^{\zeta }\right) \right] ,  \label{key3}
\end{eqnarray}%
with $\varpi \left( \epsilon \right) \rightarrow 0$ as $\epsilon \rightarrow
0$. Note that, for each $\zeta >0$, $\varpi \left( \kappa ^{\zeta }\right) $
is square integrable, we set
\begin{equation*}
\varpi _{0}\left( \zeta \right) =\mathbb{E}\left[ \varpi \left( \kappa
^{\zeta }\right) ^{2}+\varpi \left( \kappa ^{\zeta }\right) \right] .
\end{equation*}%
Hence,
\begin{equation}
\left\vert Y_{\bar{t}}^{1}-Y_{\bar{t}}^{2}\right\vert \leq C\zeta \varpi
_{0}\left( \zeta \right)  \label{key4}
\end{equation}

From the monotonicity of $G\left[ \cdot \right] ,$
\begin{eqnarray*}
\varphi \left( \bar{t},\bar{x}\right) &=&V^{-}\left( \bar{t},\bar{x}\right)
\geq G_{\bar{t},\bar{t}+\zeta }^{\bar{t},\bar{x}}\left[ V^{-}\left( \bar{t}%
+\zeta ,X_{\bar{t}+\zeta }^{\bar{t},\bar{x}}\right) \right] -\varepsilon . \\
&\geq &G_{\bar{t},\bar{t}+\zeta }^{\bar{t},\bar{x}}\left[ \varphi \left(
\bar{t}+\zeta ,X_{\bar{t}+\zeta }^{\bar{t},\bar{x}}\right) \right]
-\varepsilon .
\end{eqnarray*}%
From (\ref{key2}) and letting $\varepsilon \rightarrow 0,$
\begin{equation*}
0\geq G_{\bar{t},\bar{t}+\zeta }^{\bar{t},\bar{x}}\left[ \varphi \left( \bar{%
t}+\zeta ,X_{\bar{t}+\zeta }^{\bar{t},\bar{x}}\right) \right] -\varphi
\left( \bar{t},\bar{x}\right) =Y_{\bar{t}}^{1}.
\end{equation*}%
By (\ref{key4}) we further have $Y_{0}\left( \bar{t}\right) =Y_{\bar{t}%
}^{2}\leq C\zeta \varpi _{0}\left( \zeta \right) .$ Therefore, it follows
easily that $F\left( \bar{t},\bar{x},0,0\right) \leq 0$ and from the
definition of $F$ we see that $V^{-}$ is a viscosity supersolution of (\ref%
{HJBI}). The proof is similar for the viscosity sub-solution. \hfill $\Box $

Next, we shall prove that the HJBI equation (\ref{HJBI}) has a unique
viscosity solution. Consequently, the lower and upper value functions
coincide, since they are both viscosity solutions to (\ref{HJBI}). Thus, the
stochastic differential game admits a value.

Before introducing the comparison principle, we need the following two
technical lemmas, mainly taken from \cite{co}.

\begin{lemma}
\label{te1}Assume that \emph{(A3)} is in force. Let $U$, $V:\left[ 0,T\right]
\times \mathbb{R}^{n}\rightarrow \mathbb{R}$ a viscosity supersolution and a
viscosity subsolution to the HJBI equation (\ref{HJBI}), respectively. Let $%
\left( \hat{t},x_{0}\right) \in \left[ 0,T\right] \times \mathbb{R}^{n}$ be
such that
\begin{equation*}
V\left( \hat{t},x_{0}\right) \leq \mathcal{H}_{\sup }^{c}V\left( \hat{t}%
,x_{0}\right) ,\text{ }U\left( \hat{t},x_{0}\right) \leq \mathcal{H}_{\inf
}^{\chi }U\left( \hat{t},x_{0}\right) ,
\end{equation*}%
or%
\begin{equation*}
U\left( \hat{t},x_{0}\right) \geq \mathcal{H}_{\inf }^{\chi }U\left( \hat{t}%
,x_{0}\right) .
\end{equation*}%
Then for every $\varepsilon >0,$ there exists $\hat{x}\in \mathbb{R}^{n}$
such that
\begin{equation*}
V\left( \hat{t},x_{0}\right) -U\left( \hat{t},x_{0}\right) \leq V\left( \hat{%
t},\hat{x}\right) -U\left( \hat{t},\hat{x}\right) +\varepsilon
\end{equation*}%
and%
\begin{equation*}
V\left( \hat{t},\hat{x}\right) >\mathcal{H}_{\sup }^{c}V\left( \hat{t},\hat{x%
}\right) ,\text{ }U\left( \hat{t},\hat{x}\right) <\mathcal{H}_{\inf }^{\chi
}U\left( \hat{t},\hat{x}\right) .
\end{equation*}
\end{lemma}

\begin{lemma}
\label{te2}Assume that \emph{(A3)} is in force. Let $U$, $V:\left[ 0,T\right]
\times \mathbb{R}^{n}\rightarrow \mathbb{R}$ a viscosity supersolution and a
viscosity subsolution to the HJBI equation (\ref{HJBI}), respectively. Let $%
\left( \hat{t},\hat{x}\right) \in \left[ 0,T\right] \times \mathbb{R}^{n}$
be such that
\begin{equation*}
V\left( \hat{t},\hat{x}\right) >\mathcal{H}_{\sup }^{c}V\left( \hat{t},\hat{x%
}\right) ,\text{ }U\left( \hat{t},\hat{x}\right) <\mathcal{H}_{\inf }^{\chi
}U\left( \hat{t},\hat{x}\right) ,
\end{equation*}%
then there exists $\epsilon >0$ for which
\begin{equation*}
V\left( t,x\right) >\mathcal{H}_{\sup }^{c}V\left( t,x\right) ,\text{ }%
U\left( t,x\right) <\mathcal{H}_{\inf }^{\chi }U\left( t,x\right) ,
\end{equation*}%
where $\left( t,x\right) \in \left[ \left( \hat{t}-\delta \right) \vee
0,\left( \hat{t}+\delta \right) \wedge T\right] \times \bar{B}_{\delta
}\left( \hat{x}\right) .$
\end{lemma}

In order to get the uniqueness, we add the following assumption:

\begin{description}
\item[(A4)] Assume that $f$ is strictly monotone in $y$, that is, $f\left(
t,x,y_{1},z\right) <f\left( t,x,y_{2},z\right) ,$ for $\forall y_{1},$ $%
y_{2}\in \mathbb{R}$ with $y_{1}<y_{2},$ $\forall \left( t,x,z\right) \in %
\left[ 0,T\right] \times \mathbb{R}^{n}\times \mathbb{R}^{d}.$
\end{description}

\begin{theorem}
Let $U$, $V:\left[ 0,T\right] \times \mathbb{R}^{n}\rightarrow \mathbb{R}$ a
viscosity supersolution and a viscosity subsolution to the HJBI equation
\emph{(\ref{HJBI})}, respectively. Assume that \emph{(A1)-(A4)} are in force
and that $U$, $V$ are uniformly continuous on $\left[ 0,T\right) \times
\mathbb{R}^{n}.$ Then, we have $U\geq V$ on $\left[ 0,T\right] \times
\mathbb{R}^{n}.$
\end{theorem}

\paragraph{Proof}

We prove our result by contradiction. Suppose that
\begin{equation*}
\sup_{\left[ 0,T\right] \times \mathbb{R}^{n}}\left( V-U\right) >0.
\end{equation*}%
Fix $\theta >C_{f}>0$ where $C_{f}$ denotes the Lipschitz constant of $f.$

Define
\begin{equation*}
\bar{U}\left( t,x\right) =e^{\theta t}U\left( t,x\right) ,\text{ }\bar{V}%
\left( t,x\right) =e^{\theta t}V\left( t,x\right) ,\text{ }\left( t,x\right)
\in \left[ 0,T\right] \times \mathbb{R}^{n}.
\end{equation*}%
It is fairly easy to check that $\bar{U}\left( t,x\right) $ ($\bar{V}\left(
t,x\right) $) is a viscosity supersolusion (subsolution) to the following
HJBI equation:%
\begin{equation}
\left\{
\begin{array}{l}
\max \left\{ W-\mathcal{\bar{H}}_{\sup }^{\chi }W,\min \left\{ \theta W-%
\frac{\partial W}{\partial t}-\mathcal{L}W-\bar{f},W-\mathcal{\bar{H}}_{\inf
}^{c}W\right\} \right\} =0, \\
W\left( T,x\right) =\bar{\Phi}\left( x\right) ,\text{ }\left( t,x\right) \in %
\left[ 0,T\right) \times \mathbb{R}^{n},%
\end{array}%
\right.  \label{HJBI3}
\end{equation}%
where
\begin{eqnarray}
\mathcal{L}W\left( t,x\right) &=&\left\langle b\left( t,x\right) ,DW\left(
t,x\right) \right\rangle +\frac{1}{2}\text{tr}\left[ \sigma \sigma ^{\top
}\left( t,x\right) D^{2}W\left( t,x\right) \right] ,  \label{oper} \\
\bar{f}\left( t,x,W,DW\cdot \sigma \left( t,x\right) \right) &=&e^{\theta
t}f\left( t,x,e^{-\theta t}W,e^{-\theta t}DW\cdot \sigma \left( t,x\right)
\right) ,  \notag \\
\bar{\Phi}\left( x\right) &=&e^{\theta T}\Phi \left( x\right) ,  \notag \\
\mathcal{\bar{H}}_{\sup }^{\chi }W\left( t,x\right) &=&\sup_{z\in \mathcal{V}%
}\left[ W\left( t,x+z\right) +e^{\theta t}\chi \left( t,z\right) \right] ,
\notag \\
\mathcal{\bar{H}}_{\inf }^{c}W\left( t,x\right) &=&\inf_{y\in \mathcal{U}}%
\left[ W\left( t,x+y\right) -e^{\theta t}c\left( t,z\right) \right] .  \notag
\end{eqnarray}%
Assume that there exists $x_{0}\in \mathbb{R}^{n}$ such that $\left( \bar{U}-%
\bar{V}\right) \left( T,x_{0}\right) <0$. Then, from Lemma \ref{te1}, There
exists $\tilde{x}\in \mathbb{R}^{n}$ such that $\left( \bar{U}-\bar{V}%
\right) \left( T,\tilde{x}\right) <0.$ On the other hand, from the
subsolution property of $\bar{V}$ , we know $\bar{V}\left( T,\tilde{x}%
\right) \leq \bar{\Phi}\left( \tilde{x}\right) $. Similarly, uitilizing the
supersolution property of $\bar{U}$ , we have $\bar{U}\left( T,\tilde{x}%
\right) \geq \bar{\Phi}\left( \tilde{x}\right) $. Therefore, $\bar{U}\left(
T,\tilde{x}\right) \geq $ $\bar{V}\left( T,\tilde{x}\right) $ which leads a
contradition to $\left( \bar{U}-\bar{V}\right) \left( T,\tilde{x}\right) <0.$

Now postulate that there exists $\left( \bar{t},\bar{x}\right) \in \left[ 0,T%
\right] \times \mathbb{R}^{n}$ such that $\left( \bar{U}-\bar{V}\right)
\left( \bar{t},\bar{x}\right) <0$. Then, from Lemma \ref{te2}, there exist $%
\left( \hat{t},\hat{x}\right) \in \left[ 0,T\right] \times \mathbb{R}^{n}$
and $\delta >0$ such that%
\begin{equation*}
\sup_{I\times \bar{B}_{\delta }\left( \hat{x}\right) }\left( \bar{V}-\bar{U}%
\right) \left( t,x\right) >0
\end{equation*}%
and
\begin{equation*}
\bar{V}\left( t,x\right) >\mathcal{\bar{H}}_{\sup }^{c}V\left( t,x\right) ,%
\text{ }\bar{U}\left( t,x\right) <\mathcal{\bar{H}}_{\inf }^{\chi }\bar{U}%
\left( t,x\right) ,\text{ }\left( t,x\right) \in I\times \bar{B}_{\delta
}\left( \hat{x}\right)
\end{equation*}%
with $I:=\left[ \left( \hat{t}-\delta \right) \vee 0,\left( \hat{t}+\delta
\right) \wedge T\right] .$

We define
\begin{equation*}
\tilde{V}\left( t,x\right) =\bar{V}\left( t,x\right) -\frac{16\left\vert x-%
\hat{x}\right\vert ^{4}M}{15\delta ^{4}}\mathbf{1}_{\left\{ \left\vert x-%
\hat{x}\right\vert >\frac{\delta }{2}\right\} }+\frac{M}{15},
\end{equation*}%
where $M:=\sup_{I\times \bar{B}_{\delta }\left( \hat{x}\right) }\left( \bar{V%
}-\bar{U}\right) \left( t,x\right) .$ It is easy to check that
\begin{eqnarray*}
\left( \tilde{V}-\bar{U}\right) \left( t,x\right) &=&\left( \bar{V}-\bar{U}%
\right) \left( t,x\right) -\frac{16\left\vert x-\hat{x}\right\vert ^{4}m}{%
15\delta ^{4}}\mathbf{1}_{\left\{ \left\vert x-\hat{x}\right\vert >\frac{%
\delta }{2}\right\} }-\frac{m}{15} \\
&\leq &M.
\end{eqnarray*}%
So withoutloss of generality, we may assume that
\begin{equation*}
\left( \tilde{V}-\bar{U}\right) \left( t,x\right) \leq 0,\text{ }\left(
t,x\right) \in I\times \partial \bar{B}_{\delta }\left( \hat{x}\right) .
\end{equation*}%
Note that $\mathcal{P}^{2,+}\bar{V}\left( t,x\right) =\mathcal{P}^{2,+}%
\tilde{V}\left( t,x\right) $ for all $\left( t,x\right) \in \left[ 0,T\right]
\times \mathbb{R}^{n}$, $\tilde{V}$ can be replaced with $\bar{V}$.

Now choose $\left( t^{\prime },x^{\prime }\right) \in I\times \bar{B}%
_{\delta }\left( \hat{x}\right) $ such that
\begin{equation}
\sup_{I\times \bar{B}_{\delta }\left( \hat{x}\right) }\left( \bar{V}-\bar{U}%
\right) \left( t,x\right) =\left( \bar{V}-\bar{U}\right) \left( t^{\prime
},x^{\prime }\right) >0.  \label{contr}
\end{equation}%
Define the following text function:%
\begin{equation*}
\phi _{n}\left( t,x,y\right) =\bar{V}\left( t,x\right) -\bar{U}\left(
t,y\right) -\psi _{n}\left( t,x,y\right) ,\text{ }n\in N
\end{equation*}%
with
\begin{equation*}
\psi _{n}\left( t,x,y\right) =\frac{n}{2}\left\vert x-y\right\vert
^{2}+\left\vert x-x^{\prime }\right\vert ^{2}+\left\vert t-t^{\prime
}\right\vert ^{2}
\end{equation*}%
for every $\left( t,x,y\right) \in \left[ 0,T\right] \times \mathbb{R}%
^{n}\times \mathbb{R}^{n}$. Clearly, given any $n\geq 1,$ there exists $%
\left( t_{n},x_{n},y_{n}\right) \in I\times \bar{B}_{\delta }\left( \hat{x}%
\right) \times \bar{B}_{\delta }\left( \hat{x}\right) $ attaining the
maximum of $\phi _{n}$ on $I\times \bar{B}_{\delta }\left( \hat{x}\right)
\times \bar{B}_{\delta }\left( \hat{x}\right) $. Up to a subsequence, $%
\left( t_{n},x_{n},y_{n}\right) \in I\times \bar{B}_{\delta }\left( \hat{x}%
\right) \times \bar{B}_{\delta }\left( \hat{x}\right) \rightarrow \left(
t_{0},x_{0},y_{0}\right) \in I\times \bar{B}_{\delta }\left( \hat{x}\right)
\times \bar{B}_{\delta }\left( \hat{x}\right) $ as $n\rightarrow \infty .$
Nonetheless, for every $n\geq 1,$ we have%
\begin{equation*}
\left( \bar{V}-\bar{U}\right) \left( t^{\prime },x^{\prime }\right) =\phi
_{n}\left( t^{\prime },x^{\prime },x^{\prime }\right) \geq \phi _{n}\left(
t_{n},x_{n},y_{n}\right) .
\end{equation*}%
It yields that
\begin{eqnarray}
\left( \bar{V}-\bar{U}\right) \left( t^{\prime },x^{\prime }\right) &\leq &%
\underset{n\rightarrow \infty }{\sup \lim }\phi _{n}\left(
t_{n},x_{n},y_{n}\right)  \notag \\
&\leq &\bar{V}\left( t_{0},x_{0}\right) -\bar{U}\left( t_{0},y_{0}\right) -%
\underset{n\rightarrow \infty }{\inf \lim }n\left\vert
x_{n}-y_{n}\right\vert ^{2}  \notag \\
&&-\left\vert x_{0}-x^{\prime }\right\vert ^{2}-\left\vert t_{0}-t^{\prime
}\right\vert ^{2},  \label{opti}
\end{eqnarray}%
from which, up to a subsequence, $\underset{n\rightarrow \infty }{\inf \lim }%
n\left\vert x_{n}-y_{n}\right\vert ^{2}<\infty $. Then it follows that $%
x_{0}=y_{0}.$ From (\ref{opti}), we derive that
\begin{equation}
\left\{
\begin{array}{l}
\text{1) }\left( t_{n},x_{n},y_{n}\right) \rightarrow \left( t^{\prime
},x^{\prime },x^{\prime }\right) , \\
\text{2) }n\left\vert x_{n}-y_{n}\right\vert ^{2}\rightarrow 0 \\
\text{3) }\bar{V}\left( t_{n},x_{n}\right) -\bar{U}\left( t_{n},y_{n}\right)
\rightarrow \bar{V}\left( t^{\prime },x^{\prime }\right) -\bar{U}\left(
t^{\prime },x^{\prime }\right) ,%
\end{array}%
\right.  \label{lim}
\end{equation}%
as $n\rightarrow \infty .$ By virtue of Ishii's lemma (Theorem 8.3 in \cite%
{CIL}), up to a subsequence, we may find sequence $\left(
p_{n}^{1},q_{n}^{1},Q_{n}^{1}\right) \in \mathcal{P}^{2,+}\bar{V}\left(
t_{n},x_{n}\right) $ and $\left( p_{n}^{2},q_{n}^{2},Q_{n}^{2}\right) \in
\mathcal{P}^{2,-}\bar{U}\left( t_{n},y_{n}\right) $ such that
\begin{eqnarray*}
p_{n}^{1}-p_{n}^{2} &=&2\left( t_{n}-t^{\prime }\right) , \\
q_{n}^{1} &=&D_{x}\psi _{n}\left( t_{n},x_{n},y_{n}\right) =n\left(
x_{n}-y_{n}\right) , \\
q_{n}^{2} &=&-D_{y}\psi _{n}\left( t_{n},x_{n},y_{n}\right) =n\left(
x_{n}-y_{n}\right)
\end{eqnarray*}%
and
\begin{equation*}
\left(
\begin{array}{cc}
Q_{n}^{1} & O \\
O & -Q_{n}^{2}%
\end{array}%
\right) \leq A_{n}+\frac{1}{2n}A_{n}^{2}
\end{equation*}%
where
\begin{equation*}
A_{n}=D_{xy}\psi _{n}\left( t_{n},x_{n},y_{n}\right) =n\left(
\begin{array}{cc}
I & -I \\
-I & I%
\end{array}%
\right) .
\end{equation*}%
Then,
\begin{equation}
\left(
\begin{array}{cc}
Q_{n}^{1} & O \\
O & -Q_{n}^{2}%
\end{array}%
\right) \leq 2n\left(
\begin{array}{cc}
I & -I \\
-I & I%
\end{array}%
\right) .  \label{segma}
\end{equation}%
From $\bar{U}\left( t,x\right) $ ($\bar{V}\left( t,x\right) $) is a
viscosity supersolusion (subsolution) to the following HJBI equation (\ref%
{HJBI3}), we have%
\begin{eqnarray}
\theta \bar{U}\left( t_{n},y_{n}\right) -p_{n}^{2}-\mathcal{L}\bar{U}\left(
t_{n},y_{n}\right) -\bar{f}\left( t_{n},y_{n},e^{-\theta t_{n}}\bar{U}%
,e^{-\theta t_{n}}D\bar{U}\cdot \sigma \left( t_{n},y_{n}\right) \right)
&\geq &0,  \label{super} \\
\theta \bar{V}\left( t_{n},x_{n}\right) -p_{n}^{1}-\mathcal{L}\bar{V}\left(
t_{n},x_{n}\right) -\bar{f}\left( t_{n},x_{n},e^{-\theta t_{n}}\bar{V}%
,e^{-\theta t_{n}}D\bar{V}\cdot \sigma \left( t,x_{n}\right) \right) &\leq
&0,  \label{sub}
\end{eqnarray}%
where $\mathcal{L}$ is defined in (\ref{oper}). From (\ref{super}) and (\ref%
{sub}), we immediately get%
\begin{eqnarray}
&&\theta \bar{V}\left( t_{n},x_{n}\right) -\theta \bar{U}\left(
t_{n},x_{n}\right) +p_{n}^{2}-p_{n}^{1}+\mathcal{L}\bar{U}\left(
t_{n},x_{n}\right) -\mathcal{L}\bar{V}\left( t_{n},y_{n}\right)  \notag \\
&&+\bar{f}\left( t_{n},y_{n},e^{-\theta t_{n}}\bar{U}\left(
t_{n},y_{n}\right) ,e^{-\theta t_{n}}q_{n}^{2}\cdot \sigma \left(
t_{n},y_{n}\right) \right)  \notag \\
&&-\bar{f}\left( t_{n},x_{n},e^{-\theta t_{n}}\bar{V}\left(
t_{n},x_{n}\right) ,e^{-\theta t_{n}}q_{n}^{1}\cdot \sigma \left(
t_{n},x_{n}\right) \right)  \notag \\
&\leq &0.  \label{u0}
\end{eqnarray}%
Clearly,
\begin{equation}
p_{n}^{1}-p_{n}^{2}\rightarrow 0,\text{ as }n\rightarrow \infty .  \label{u1}
\end{equation}%
and
\begin{equation}
\left\langle b\left( t_{n},y_{n}\right) ,n\left( x_{n}-y_{n}\right)
\right\rangle -\left\langle b\left( t_{n},x_{n}\right) ,n\left(
x_{n}-y_{n}\right) \right\rangle \rightarrow 0,  \label{u2}
\end{equation}%
since (1)-(2) in (\ref{lim}).

For simplicity, set $\sigma _{1}=\sigma \left( t_{n},y_{n}\right) ,$ $\sigma
_{2}=\sigma \left( t_{n},x_{n}\right) .$ We deal with%
\begin{eqnarray}
&&\frac{1}{2}\text{\textrm{Tr}}\left( \sigma \sigma ^{\ast }\left(
t_{n},y_{n}\right) Q_{n}^{1}\right) -\frac{1}{2}\text{\textrm{Tr}}\left(
\sigma \sigma ^{\ast }\left( t_{n},x_{n}\right) Q_{n}^{2}\right)  \notag \\
&=&\frac{1}{2}\text{\textrm{Tr}}\left(
\begin{array}{cc}
\sigma _{1}\sigma _{1}^{\top } & \sigma _{1}\sigma _{2}^{\top } \\
\sigma _{2}\sigma _{1}^{\top } & \sigma _{2}\sigma _{2}^{\top }%
\end{array}%
\right) \left(
\begin{array}{cc}
Q_{n}^{1} & 0 \\
0 & -Q_{n}^{2}%
\end{array}%
\right)  \notag \\
&\leq &n\text{\textrm{Tr}}\left(
\begin{array}{cc}
\sigma _{1}\sigma _{1}^{\top } & \sigma _{1}\sigma _{2}^{\top } \\
\sigma _{2}\sigma _{1}^{\top } & \sigma _{2}\sigma _{2}^{\top }%
\end{array}%
\right) \left(
\begin{array}{cc}
I & -I \\
-I & I%
\end{array}%
\right)  \notag \\
&\leq &n\text{\textrm{Tr}}\left[ \sigma _{1}\sigma _{1}^{\top }-\sigma
_{1}\sigma _{2}^{\top }-\sigma _{2}\sigma _{1}^{\top }+\sigma _{2}\sigma
_{2}^{\top }\right]  \notag \\
&=&n\text{\textrm{Tr}}\left[ \left( \sigma _{1}-\sigma _{2}\right) \left(
\sigma _{1}-\sigma _{2}\right) ^{\top }\right]  \notag \\
&\leq &nC\left\vert \sigma _{1}-\sigma _{2}\right\vert ^{2}  \notag \\
&\leq &nC\left\vert x_{n}-y_{n}\right\vert ^{2}\rightarrow 0,\text{ as }%
n\rightarrow \infty ,  \label{u3}
\end{eqnarray}%
where we have used the the assumption that Lipschitz condition on $\sigma $,
(\ref{segma}) and (2) in (\ref{lim}).

Thanks to (3) in (\ref{lim}), the left-hand side of inequality (\ref{u0})
goes to $\theta \left[ \bar{V}\left( t^{\prime },x^{\prime }\right) -\bar{U}%
\left( t^{\prime },x^{\prime }\right) \right] $, as $n\rightarrow \infty ,$
moreover, by (A4), we have
\begin{eqnarray*}
\theta \left[ \bar{V}\left( t^{\prime },x^{\prime }\right) -\bar{U}\left(
t^{\prime },x^{\prime }\right) \right] &\leq &\bar{f}\left( t^{\prime
},x^{\prime },e^{-\theta t^{\prime }}\bar{U}\left( t^{\prime },x^{\prime
}\right) ,0\right) -\bar{f}\left( t^{\prime },x^{\prime },e^{-\theta
t^{\prime }}\bar{V}\left( t^{\prime },x^{\prime }\right) ,0\right) \\
&<&0
\end{eqnarray*}%
which leads to a contradiction to (\ref{contr}). Our proof is thus
completed.\hfill $\Box $

\begin{remark}
\label{rk}To get a uniqueness result for viscosity solution of (\ref{HJBI}),
we adapt some techniques from \cite{co}. We have to mention that there is
another approach developed by Barles, Buckdahn and Pardoux \cite{BBP}. The
value function can be considered in given class of continuous functions
satisfying
\begin{equation*}
\lim_{\left\vert x\right\vert \rightarrow \infty }\left\vert u\left(
t,x\right) \right\vert \exp \left\{ -A\left[ \log \left( \left\vert
x\right\vert \right) \right] ^{2}\right\} =0,
\end{equation*}%
uniformly for $t\in \left[ 0,T\right] ,$ for some $A>0.$ The space of
continuous functions endowed with a growth condition is slightly weaker than
the assumption of polynomial growth but more restrictive than that of
exponential growth. This growth condition was first introduced by Barles,
Buckdahn, and Pardoux \cite{BBP} to prove the uniqueness of the viscosity
solution of an integro-partial differential equation associated with a
decoupled FBSDEs with jumps. It has been shown in \cite{BBP} that this kind
of growth condition is optimal for the uniqueness and can, in general, not
be weakened. These techniques have been applied in \cite{BJ, WY} for the
uniqueness for viscosity solutions of recursive control of the obstacle
constraint problem and Hamilton-Jacobi-Bellman-Isaacs equations related to
stochastic differential games, respectively. However, as you may have
observed, in our HJBI equation, there appears two obstacles, which are
implicit obstacles, in the sense that they depend on $V^{-}$. It is worth to
pointing out that the smooth supersolution built in \cite{BBP}, namely
\begin{equation*}
\chi \left( t,x\right) =\exp \left[ \left( \check{C}\left( T-t\right)
+A\right) \psi \left( x\right) \right] ,
\end{equation*}%
whilst
\begin{equation*}
\psi \left( x\right) =\left[ \log \left( \left\vert x\right\vert
^{2}+1\right) ^{\frac{1}{2}}+1\right] ^{2},
\end{equation*}%
where $\check{C}$ and $A$ are positive constants. One can show
\begin{equation*}
\min_{v\in U}\left\{ \mathcal{L}\left( t,x,v\right) \chi \left( t,x\right)
+C\left\vert \chi \right\vert +C\left\vert \nabla \chi \sigma \left(
t,x,v\right) \right\vert \right\} \leq 0,
\end{equation*}%
where $C$ is the Lipschitz constant of $f.$ Following the idea in \cite{BBP}%
, whenever considering the difference of $u_{1}-u_{2}$ where $u_{1}$ ($u_{2}$%
) is a subsolution (supersolution) of (\ref{HJBI}). It is hard to check the
obstacles of viscosity solution (\ref{vis1})-(\ref{TT2}). This is the reason
we borrow the idea from Fleming, Soner \cite{FS} and Cosso \cite{co} to
handle the uniqueness.
\end{remark}

\begin{remark}
As observed in our paper, we put somewhat strong assumptions on
coefficients, namely, \textit{boundedness}. On the one hand, it simplifies
our proof of existence. Recently, El Asri and Mazid \cite{AM} also
investigate the solution to the zero-sum stochastic differential games, but
under rather weak assumptions on the cost functions ($c$ and $\chi $ are not
decreasing in time). In the future, we shall adopt the idea developed by El
Asri and Mazid \cite{AM} to exploit the recursive utilities.
\end{remark}

\section{Concluding remarks}

\label{sec5}

In this paper, we study on zero-sum stochastic differential games in the
framework of backward stochastic differential equations on a finite time
horizon with both players adopting impulse controls. By means of stochastic
backward semigroups and comparison theorem of BSDE, we prove a dynamic
programming principle for both the upper and the lower value functions of
the game. The upper and the lower value functions are then shown to be the
unique viscosity solutions of the Hamilton-Jacobi-Bellman-Isaacs equations
with a double-obstacle. As a result, the uniqueness implies that the upper
and lower value functions coincide and the game admits a value. In future
work, we plan to relax our assumptions and to try to find a smooth solution
for the HJBI (\ref{HJBI}) in order to obtain uniqueness as Remark \ref{rk}.
Besides, as in Zhang \cite{ZF}, we will consider problems in which one
player adopts impulse controls and the other adopts continuous controls,
finite/infinite horizons, \emph{etc}. These possible
extensions promise
to be interesting research directions. We shall response these challenging
topics in our future work.

\appendix

\section{The Proof of Lemma \protect\ref{deter}}

\label{APP}

\paragraph{Proof.}

We adopt the idea from \cite{BuckLi}. Let $H$ denote the Cameron--Martin
space of all absolutely continuous elements $h\in \Omega $ whose derivative $%
\dot{h}$ belongs to $L^{2}\left( \left[ 0,T\right] ,\mathbb{R}^{d}\right) $.
For any $h\in H$, we define the mapping $\tau _{h}\omega :=\omega +h$, $%
\omega \in \Omega $. Clearly, $\tau _{h}:\Omega \rightarrow \Omega $ is a
bijection, and its law is given by
\begin{equation*}
P\circ \left( \tau _{h}\right) ^{-1}=\exp \left\{ \int_{0}^{T}\dot{h}_{s}%
\mathrm{d}W_{s}-\frac{1}{2}\int_{0}^{T}\left\vert \dot{h}_{s}\right\vert ^{2}%
\mathrm{d}s\right\} P.
\end{equation*}
Let $\left( t,x\right) \in \left[ 0,T\right] \times \mathbb{R}^{n}$ be
arbitrarily fixed, and put $H_{t}=\left\{ h\in H|h\left( \cdot \right)
=h\left( \cdot \wedge t\right) \right\} $. Let $u\in \mathcal{U}_{t,T},$ $%
v\in \mathcal{V}_{t,T}$, and $h\in H_{t}$, we first show that $J\left( t,x;u,v\right) \left(
\tau _{h}\right) =J\left( t,x;u\left( \tau _{h}\right) ,v\left( \tau
_{h}\right) \right) $, $P$-a.s. Indeed, substitute the transformed control
processes $u\left( \tau _{h}\right) $ and $v\left( \tau _{h}\right) $ for $u$
and $v$ into FBSDEs (\ref{SDE1})-(\ref{BSDE1}) and take the
Girsanov transformation to (\ref{SDE1})-(\ref{BSDE1}), finally compare the obtained equation with the previous ones.
Then from the uniqueness of the solution of (\ref{SDE1})-(\ref{BSDE1}), we conclude with
\begin{eqnarray*}
X_{s}^{t,x;u,v}\left( \tau _{h}\right)  &=&X_{s}^{t,x;u\left( \tau
_{h}\right) ,v\left( \tau _{h}\right) }, \\
Y_{s}^{t,x;u,v}\left( \tau _{h}\right)  &=&Y_{s}^{t,x;u\left( \tau
_{h}\right) ,v\left( \tau _{h}\right) }, \\
Z_{s}^{t,x;u,v}\left( \tau _{h}\right)  &=&Z_{s}^{t,x;u\left( \tau
_{h}\right) ,v\left( \tau _{h}\right) }
\end{eqnarray*}%
for any $s\in \left[ t,T\right] $, $P$-a.s., which indicates that $J\left(
t,x;u,v\right) \left( \tau _{h}\right) =J\left( t,x;u\left( \tau _{h}\right)
,v\left( \tau _{h}\right) \right) ,$ $P$-a.s. For $\beta \in \mathcal{B}%
_{t,T},$ $h\in H_{t}$, let $\beta ^{h}\left( u\right) :=\beta \left( u\left(
\tau _{-h}\right) \right) \left( \tau _{h}\right) $, $u\in \mathcal{U}_{t,T}$%
. Then $\beta ^{h}\in \mathcal{B}_{t,T}$, which makes $\mathcal{U}_{t,T}$
into $\mathcal{V}_{t,T}.$ Moreover, it is easy to check that this mapping is
nonanticipating and verify
\begin{equation*}
\left\{ \sup_{u\in \mathcal{U}_{t,T}}J\left( t,x;u,\beta \left( u\right)
\right) \right\} \left( \tau _{h}\right) =\sup_{u\in \mathcal{U}%
_{t,T}}\left\{ J\left( t,x;u,\beta \left( u\right) \right) \left( \tau
_{h}\right) \right\} ,\text{ }P\text{-a.s.}
\end{equation*}%
Now let $h\in H_{t},$
\begin{eqnarray*}
V^{-}\left( t,x\right) \left( \tau _{h}\right) &=&\inf_{\beta \in \mathcal{B}%
_{t,T}}\sup_{u\in \mathcal{U}_{t,T}}\left\{ J\left( t,x;u,\beta \left(
u\right) \right) \left( \tau _{h}\right) \right\} \\
&=&\inf_{\beta \in \mathcal{B}_{t,T}}\sup_{u\in \mathcal{U}_{t,T}}\left\{
J\left( t,x;u\left( \tau _{h}\right) ,\beta ^{h}\left( u\left( \tau
_{h}\right) \right) \right) \right\} \\
&=&\inf_{\beta \in \mathcal{B}_{t,T}}\sup_{u\in \mathcal{U}_{t,T}}\left\{
J\left( t,x;u,\beta ^{h}\left( u\right) \right) \right\} \\
&=&\inf_{\beta \in \mathcal{B}_{t,T}}\sup_{u\in \mathcal{U}_{t,T}}\left\{
J\left( t,x;u,\beta \left( u\right) \right) \right\} \\
&=&V^{-}\left( t,x\right) ,\text{ }P\text{-a.s.,}
\end{eqnarray*}%
which holds even for all $h\in H$. Recall the definition of the filtration,
the $\mathcal{F}_{t}$-measurable random variable $V^{-}\left( t,x\right)
\left( \omega \right) $, $\omega \in \Omega $, depends only on the
restriction of $\omega $ to the time interval $\left[ 0,t\right] $. We
complete our proof with help of Lemma 3.4 in \cite{BuckLi}.\hfill $\Box $

\noindent \textbf{Acknowledgements.} The author wishes to thank the AE and
referees for their careful reading and helpful comments that improved the
first version of the paper. The author also thanks the department of applied
mathematics, The Hong Kong Polytechnic University, for its hospitality
during his visit in Jan. 2019.

\end{document}